\documentclass[12pt]{article}
\usepackage[margin=1in, a4paper]{geometry}
\usepackage{amsmath,enumerate,amsfonts,color,amssymb,amsthm, verbatim}

\usepackage{multirow}
\usepackage{hyperref}
\usepackage[normalem]{ulem}

\usepackage{todonotes}

\usepackage{slashed}

\usepackage{tikz}

\def\RR{{\mathbb R}}

\def\eps{{\varepsilon}}

\numberwithin{equation}{section}

\newtheorem{theorem} {\sc  Theorem\rm}[section]

\newtheorem{corollary} [theorem] {\sc  Corollary\rm}
\newtheorem{lemma} [theorem] {\sc  Lemma\rm}
\newtheorem{proposition} [theorem] {\sc  Proposition\rm}

\newcounter{marnote}

\DeclareFontFamily{OT1}{rsfs}{}
\DeclareFontShape{OT1}{rsfs}{m}{n}{ <-7> rsfs5 <7-10> rsfs7 <10-> rsfs10}{}
\DeclareMathAlphabet{\mycal}{OT1}{rsfs}{m}{n}

\def\be{\begin{equation}}
\def\ee{\end{equation}}

\newcommand{\R}{\mathbb{R}}

\def\be{\begin{equation}}
\def\ee{\end{equation}}
\def\bea#1\eea{\begin{align}#1\end{align}}

\def\mcA{{\mycal A}}

\def\mcE{{\mycal E}}
\def\mcF{{\mycal F}}

\begin{document}
\title{Global minimality of the degree-one Ginzburg-Landau vortex solution in the unit ball in dimensions $n \geq 2$}

\author{Radu Ignat\thanks{Institut de Math\'ematiques de Toulouse, UMR 5219, Universit\'e de Toulouse, CNRS, UPS
IMT, F-31062 Toulouse Cedex 9, France. Email: Radu.Ignat@math.univ-toulouse.fr
}~, Luc Nguyen\thanks{Mathematical Institute and St Edmund Hall, University of Oxford, Andrew Wiles Building, Radcliffe Observatory Quarter, Woodstock Road, Oxford OX2 6GG, United Kingdom. Email: luc.nguyen@maths.ox.ac.uk}}

\date{}

\maketitle
\begin{abstract}
We consider the problem of minimizing the Ginzburg-Landau functional among $\RR^n$-valued maps from the unit ball $B \subset \RR^n$ with vortex boundary condition $u(x) = x$ on $\partial B$. We show that, for every Ginzburg-Landau parameter $\eps > 0$, there exists a unique minimizer given exactly by the rotationally symmetric vortex solution in dimension $n\geq 2$.


\end{abstract}

\setcounter{tocdepth}{1}
\tableofcontents

\section{Introduction}

A longstanding open problem in the Ginzburg-Landau theory asks if the degree-one vortex solution in the unit ball is minimizing and, moreover, if it is the unique global minimizer in dimension $n\geq 2$. The answer is affirmative in dimensions $n \ge 7$ (see Ignat, Nguyen, Slastikov and Zarnescu \cite{INSZ_CRAS18, INSZ-ENS}), and in dimension $n = 2$ (see the very recent work Chen, Liu, Wei and Yang \cite{CLWY-Arxiv}). In this paper, we settle the remaining case $3 \leq n \leq 6$, also in the affirmative. Our argument provides a unified proof in all dimensions $n \geq 2$.

For $\eps>0$, consider the Ginzburg-Landau energy functional
\[
\mcE_\eps(u) = \int_{B} \Big[\frac{1}{2}|\nabla u|^2 + \frac{1}{4\eps^2} (1 - |u|^2)^2\Big]\,dx
\]
where $B$ denotes the unit ball in $\RR^n$, $n \geq 2$, and $u$ belongs to
\[
\mcA = \Big\{u \in H^1(B,\RR^n): u(x) = x \text{ on } \partial B\Big\}, \quad n \geq 2.
\]
A special critical point of $\mcE_\eps$ in $\mcA$ known as the degree-one vortex solution is given by
\[
u_\eps(x) = f_\eps(|x|) \frac{x}{|x|} \quad \text{ in } B,
\]
where the scalar radial profile $f_\eps$ is the unique solution to the ODE problem in the radial variable $r=|x|$:
\begin{equation}
\begin{cases}
f_\eps'' + \frac{n-1}{r} f_\eps' - \frac{n-1}{r^2} f_\eps = - \frac{1}{\eps^2} (1 - f_\eps^2) f_\eps \quad \text{ in } (0,1),\\
f_\eps(0) = 0, f_\eps(1) = 1.
\end{cases}
	\label{Eq:feps}
\end{equation}
In particular, $0<f_\eps<1$ and $f'_\eps>0$ in $(0,1)$ (see \cite{tang, Hervex2, ODE_INSZ}).

\begin{theorem}\label{Thm:Main}
Let $n \geq 2$. Then, for every $\eps > 0$, $u_\eps$ is the unique minimizer of $\mcE_\eps$ in $\mcA$.
\end{theorem}

We briefly recall the history of the problem. The question of minimality of $u_\eps$ in the ball $B$ and also of the entire vortex solution in $\RR^n$ was raised in dimension $n = 2$ in the seminal book Bethuel, Brezis and H\'elein \cite[Problem 10]{vortices} and in higher dimension in Brezis \cite[Section 2]{Brezis99_PSPM}, \cite[Section 2]{Brezis23_Lincei}. The case of the entire vortex solution was resolved in dimension $n = 2$ by Mironescu \cite{Mironescu_symmetry}, in dimension $n = 3$ by Millot and Pisante \cite{Mil-Pis}, and in dimensions $n \geq 4$ by Pisante \cite{Pisante11-JFA}. On the unit ball, when $\eps$ is sufficiently large, $\mcE_\eps$ is strictly convex and so $u_\eps$ is the unique minimizer of $\mcE_\eps$ in $\mcA$. The minimality of $u_\eps$ for sufficiently small $\eps > 0$ in dimension $n = 2$ was proved in the monograph Pacard and Rivi\`ere \cite{PacRivBook}. The minimality of $u_\eps$ for all $\eps > 0$ was shown in dimensions $n \geq 7$ in \cite{INSZ_CRAS18, INSZ-ENS}. For gradient fields, this was lowered to $n \geq 4$ in Ignat, Nahon and Nguyen \cite{INN25}. Most recently, the two-dimensional case was settled for all $\eps > 0$ in \cite{CLWY-Arxiv}. 

For further works on the global minimality or local minimality of vortex solutions in the Ginzburg-Landau theory as well as related models, see \cite{BrezisCoronLieb, Pino-Falmer-Kow, GarPetSpa25-M3AS, GolovatyBerlyand-CVPDE02, Gustafson, HangLin01-ActaSin, IN-IHP24, INSZ3, IgnatRus24, JagerKaul83-JRAM, LamyMarconi23-SNS, LiMelcher18-JFA, LiebLoss95-JEDP, Lin-CR87, Mironescu-radial, Sandier98, SandierShafrir94-CVPDE} and the references therein. For further reading on vortices in the Ginzburg-Landau theory, see, e.g., the monographs \cite{vortices, PacRivBook, SandierSerfatyBook}.

Let us discuss some elements of our proof. Roughly speaking, our proof combines the calibration approach in \cite{INSZ_CRAS18, INSZ-ENS} with the entire-vortex comparison introduced in the two-dimensional argument of \cite{CLWY-Arxiv}. To begin with, we recall from \cite{INSZ_CRAS18, INSZ-ENS} the inequality
\begin{align}
\mcE_\eps(u) - \mcE_\eps(u_\eps) &\geq \frac{1}{2} \mcF_\eps(u - u_\eps)\nonumber\\
& := \frac{1}{2} \int_B \Big[|\nabla (u - u_\eps)|^2 - \frac{1}{\eps^2}(1 - f_\eps^2)|u - u_\eps|^2\Big]\,dx\nonumber\\
& = \frac{1}{2} \int_B f_\eps^2\Big[\Big|\nabla \Big(\frac{u - u_\eps}{f_\eps}\Big)\Big|^2 - \frac{n-1}{r^2} \frac{|u - u_\eps|^2}{f_\eps^2}\Big]\,dx
\label{Eq:Cal1}
\end{align}
for all $u \in u_\eps + C_c^\infty(B \setminus \{0\}, \RR^n)$. It was shown that, when $n \geq 7$, the quadratic functional $\mcF_\eps$ is positive definite on $H_0^1(B, \RR^n)$ and hence the vortex solution $u_\eps$ is the unique minimizer of $\mcE_\eps$ in $\mcA$. The inequality \eqref{Eq:Cal1} can therefore be viewed as a kind of calibration. However, in dimensions $2 \leq n \leq 6$, the functional $\mcF_\eps$ is no longer positive definite when $\eps$ is sufficiently small (see \cite[Lemma 2.3(c)]{IN-IHP24}), and a strengthened calibration inequality is desirable.

To construct such a calibration inequality, we make use of the minimizing property of the entire vortex solution in $\R^n$, as in \cite{CLWY-Arxiv}. For $\eps>0$, let 
\[
U_\eps(x) = F_\eps(|x|) \frac{x}{|x|} \quad \textrm{ in } \R^n,
\]
be the entire vortex solution where scalar radial profile $F_\eps$ is the unique solution to the problem
\begin{equation}
\begin{cases}
F_\eps'' + \frac{n-1}{r} F_\eps' - \frac{n-1}{r^2} F_\eps = - \frac{1}{\eps^2} (1 - F_\eps^2) F_\eps \text{ in } (0,\infty),\\
F_\eps(0) = 0, F_\eps(\infty) = 1.
\end{cases}
\label{Eq:Feps}
\end{equation}
Note that $f_\eps > F_\eps$ in $(0,1]$, $0<F_\eps<1$ and $F'_\eps>0$ in $(0, \infty)$ (see \cite{tang, Hervex2, ODE_INSZ}). The minimizing property of $U_\eps$ (see the aforementioned works \cite{Mil-Pis, Mironescu_symmetry, Pisante11-JFA}) gives 
\begin{equation}
\mcE_\eps\Big(\frac{F_\eps}{f_\eps} u\Big) - \mcE_\eps(U_\eps) \geq 0 \qquad \forall~ u \in \mcA.
	\label{Eq:Cal2}
\end{equation}
We then obtain the following analogue $\widetilde{\mcF}_\eps$ of $\mcF_\eps$:
\begin{align}
&\Delta_\eps(u):= \Big[\mcE_\eps(u) - \mcE_\eps(u_\eps)\Big] - \Big[\mcE_\eps\Big(\frac{F_\eps}{f_\eps} u\Big) - \mcE_\eps(U_\eps)\Big] \geq \frac{1}{2} \widetilde{\mcF}_\eps(u - u_\eps)\nonumber\\
&\quad  := \frac{1}{2} \int_B \Big[|\nabla (u - u_\eps)|^2 - \Big|\nabla\Big(\frac{F_\eps}{f_\eps} (u - u_\eps)\Big)\Big|^2 - \frac{f_\eps^2 - F_\eps^2 - f_\eps^4 + F_\eps^4}{\eps^2 f_\eps^2}  |u - u_\eps|^2\Big]\,dx\nonumber\\
&\quad  = \frac{1}{2} \int_B (f_\eps^2 - F_\eps^2) \Big[\Big|\nabla \Big(\frac{u - u_\eps}{f_\eps}\Big)\Big|^2 - \frac{n-1}{r^2} \frac{|u - u_\eps|^2}{f_\eps^2}\Big]\,dx 
\label{Eq:Cal3}
\end{align}
provided $u \in u_\eps + C_c^\infty(B\setminus \{0\},\RR^n)$. 
(For a related statement in dimension $n = 2$, see \cite[Proposition 2.7]{CLWY-Arxiv}.) Comparing with \eqref{Eq:Cal1}, one notes that the weight $f_\eps^2$ in \eqref{Eq:Cal1} is now replaced by $f_\eps^2 - F_\eps^2$.
 
In view of \eqref{Eq:Cal2}, the minimality of the vortex solution $u_\eps$ is established once we can show that the quadratic form $\widetilde{\mcF}_\eps$ is positive definite. Following \cite{CLWY-Arxiv}, we introduce the auxiliary function
\begin{equation}
\phi_\eps = \frac{F_\eps^2}{r^2 (f_\eps^2 - F_\eps^2)^{1/2}}>0 \quad \textrm{ in } (0, 1],
	\label{Eq:phi}
\end{equation}
and rewrite $\widetilde{\mcF}_\eps$ in the form
\begin{equation}
\widetilde{\mcF}_\eps(w) = \int_{B} (f_\eps^2 - F_\eps^2)\Big[ \phi_\eps^2 \Big|\nabla\Big(\frac{w}{f_\eps \phi_\eps}\Big)\Big|^2 + \frac{S_\eps |w|^2}{f_\eps^2}\Big]\,dx, \quad w \in H_0^1(B, \RR^n),
	\label{Eq:tmcFS}
\end{equation}
for a certain $S_\eps$ (see \eqref{Eq:S} for its definition). The main technical hurdle is then to prove the positivity of the scalar function $S_\eps$ arising in the calibration representation \eqref{Eq:tmcFS}, which immediately gives the positive definiteness of $\widetilde{\mcF}_\eps$. Our strategy to prove that $S_\eps > 0$ differs substantially from that in \cite{CLWY-Arxiv}. 
We interpolate $f_\eps$ and $F_\eps$ by the solutions $f_{\eps,a}$ of the problem
\begin{equation}
\begin{cases}
f_{\eps,a}'' + \frac{n-1}{r} f_{\eps,a}' - \frac{n-1}{r^2} f_{\eps,a} = - \frac{1}{\eps^2} (1 - f_{\eps,a}^2) f_{\eps,a} \text{ in } (0,1),\\
f_{\eps,a}(0) = 0, f_{\eps,a}(1) = a,
\end{cases}
	\label{Eq:fepsa}
\end{equation}
where $a$ varies in $[F_\eps(1),1]$. For each such $f_{\eps,a}$, we define a corresponding function $S_{\eps,a}$ (see \eqref{Eq:Sepsa}) so that $S_\eps = S_{\eps,1}$. Our proof of the positivity of $S_\eps$ is split into showing that $S_{\eps,a}$ is increasing in $a$ and that $S_{\eps,F_\eps(1)} > 0$. This has the advantage of placing fine ODE analysis largely on the entire vortex solution, where we can make use of some of the machinery developed in Ignat, Nguyen, Slastikov and Zarnescu \cite{ODE_INSZ}.

Although our proof works in all dimensions $n \geq 2$, a few algebraic estimates require separate arguments for $n = 2$, $3 \leq n \leq 7$, and $n \geq 8$, according to the sign of certain coefficients. See Lemmas \ref{Lem:psiUB} and \ref{Lem:xy>}.

The rest of the paper is organized as follows. In Section \ref{Sec2}, we give the proof of the calibration inequality \eqref{Eq:Cal3} and reduce the proof of the positive definiteness of $\widetilde{\mcF}_\eps$ to that of the positivity of $S_{\eps,a}$. In Section \ref{Sec3}, we establish refined estimates for the entire vortex solution and its linearization needed in our study of $S_{\eps,a}$. Finally, in Section \ref{Sec4}, we prove the positivity of $S_{\eps,a}$ and, in particular, of $S_\eps$.

\section{Reduction to ODE analysis}\label{Sec2}

In this section, we prove the calibration inequality \eqref{Eq:Cal3}, define $S_\eps$, establish the representation \eqref{Eq:tmcFS}, and reduce the study of the positive definiteness of $\widetilde{\mcF}_\eps$ to the study of the positivity of a family of scalar functions $S_{\eps,a}$ as mentioned in the introduction.

\begin{lemma}\label{Lem:CalS1}
Let $n \geq 2$ and $\eps > 0$. Then \eqref{Eq:Cal3} holds for $u \in u_\eps + C_c^\infty(B\setminus \{0\},\RR^n)$.
\end{lemma}

\begin{proof}
We compute as in \cite{INSZ_CRAS18}. Let $w = u - u_\eps \in  C_c^\infty(B\setminus \{0\},\RR^n)$. First, we have 
\begin{align*}
\mcE_\eps(u) - \mcE_\eps(u_\eps)
	&= \int_B \Big[\nabla u_\eps : \nabla w - \frac{1}{\eps^2} (1 - |u_\eps|^2) u_\eps \cdot w\Big]\,dx\\
		&\qquad + \int_B \Big[\frac{1}{2}|\nabla w|^2 - \frac{1}{2\eps^2} (1 - |u_\eps|^2) |w|^2 + \frac{1}{4\eps^2}(2 u_\eps \cdot w + |w|^2)^2\Big]\,dx.
\end{align*}
The first integral vanishes due to the criticality of $u_\eps$ yielding
\begin{align*}
\mcE_\eps(u) - \mcE_\eps(u_\eps)
	&= \frac{1}{2}\int_B \Big[ |\nabla w|^2 - \frac{1}{ \eps^2} (1 - f_\eps^2) |w|^2\\
		&\qquad + \frac{2}{\eps^2} f_\eps^2 (n \cdot w)^2 + \frac{2}{ \eps^2}  f_\eps |w|^2(n \cdot w ) +  \frac{1}{2\eps^2} |w|^4\Big]\,dx,
\end{align*}
where $n(x) = \frac{x}{|x|}$. Similarly,  as $\frac{F_\eps}{f_\eps}u-U_\eps=\frac{F_\eps}{f_\eps}w$,  
\begin{align*}
\mcE_\eps\Big(\frac{F_\eps}{f_\eps}u\Big) - \mcE_\eps(U_\eps)
	&= \frac{1}{2}\int_B \Big[|\nabla \Big(\frac{F_\eps}{f_\eps}w\Big)|^2 - \frac{1}{\eps^2}\frac{F_\eps^2}{f_\eps^2} (1 - F_\eps^2) |w|^2\\
		&\qquad + \frac{2}{\eps^2} \frac{F_\eps^4}{f_\eps^2} (n \cdot w)^2 + \frac{2}{ \eps^2}  \frac{F_\eps^4}{f_\eps^3}  |w|^2(n \cdot w ) +  \frac{1}{2\eps^2}\frac{F_\eps^4}{f_\eps^4} |w|^4  \Big]\,dx.
\end{align*}
Therefore,  since $f_\eps>F_\eps$ in $(0,1)$, 
\begin{align}
\nonumber
\Delta_\eps(u)  
	&=\Big[\mcE_\eps(u) - \mcE_\eps(u_\eps)\Big] - \Big[\mcE_\eps\Big(\frac{F_\eps}{f_\eps} u\Big) - \mcE_\eps(U_\eps)\Big]
\nonumber	\\ & = \frac{1}{2} \int_B \Big[|\nabla w|^2 - \Big|\nabla\Big(\frac{F_\eps}{f_\eps}w\Big)\Big|^2 - \frac{1}{\eps^2} \frac{f_\eps^2 - F_\eps^2}{f_\eps^2} (1 - f_\eps^2 - F_\eps^2)|w|^2\Big]\,dx \\\nonumber	
		&\qquad  + \frac{2}{\eps^2} \frac{1}{f_\eps^4}(f_\eps^4 - F_\eps^4) \big(f_\eps  n \cdot w  + \frac{1}{2}|w|^2\big)^2\Big]\,dx\\
\label{eq:ftilde}	&  \geq   \frac{1}{2} \int_B \Big[|\nabla w|^2 - \Big|\nabla\Big(\frac{F_\eps}{f_\eps}w\Big)\Big|^2 - \frac{1}{\eps^2} \frac{f_\eps^2 - F_\eps^2 - f_\eps^4 + F_\eps^4}{f_\eps^2}  |w|^2\Big]\,dx = \frac{1}{2} \widetilde{\mcF}_\eps(w).
\end{align}

It remains to prove for $\tilde w = \frac{w}{f_\eps}$ that
\begin{equation}
\widetilde{\mcF}_\eps(w) = \int_B (f_\eps^2 - F_\eps^2)\Big[|\nabla \tilde w|^2 - \frac{n-1}{r^2} |\tilde w|^2\Big]\,dx.
	\label{Eq:Cal4}
\end{equation}
Indeed, we compute using \eqref{Eq:feps} and the  Hardy decomposition $w=f_\eps \tilde w$ for the operator $-\Delta$ (see e.g. \cite[Lemma A.1]{INSZ3}),
\begin{align*}
\int_B |\nabla w|^2\,dx 
	&		= \int_B \Big[f_\eps^2 |\nabla \tilde w|^2 +  |\tilde w|^2 f_\eps \cdot (-\Delta f_\eps)   \Big]\,dx \\
	&= \int_B \Big[f_\eps^2 |\nabla \tilde w|^2 - \frac{n-1}{r^2} f_\eps^2|\tilde w|^2    + \frac{1}{\eps^2} f_\eps^2(1 - f_\eps^2) |\tilde w|^2  \Big]\,dx.
\end{align*}
Similarly, by \eqref{Eq:Feps},
\begin{align*}
\int_B \Big|\nabla \Big(\frac{F_\eps}{f_\eps} w\Big)|^2\,dx 
	&= \int_B |\nabla (F_\eps \tilde w)|^2\,dx\\
		 &= \int_B \Big[F_\eps^2 |\nabla \tilde w|^2 - \frac{n-1}{r^2} F_\eps^2|\tilde w|^2 + \frac{1}{\eps^2} F_\eps^2(1 - F_\eps^2) |\tilde w|^2  \Big]\,dx.
\end{align*}
Putting these two identities into the definition of $\widetilde{\mcF}_\eps(w)$ we obtain \eqref{Eq:Cal4}.
\end{proof}

Define the scalar radial functions
\begin{align}
q_\eps &= \frac{r F_\eps'}{F_\eps} \quad \textrm{ in } (0, \infty),
	\label{Eq:q}\\
S_\eps &= \frac{1}{\eps^2}(1 + f_\eps^2 - F_\eps^2) - \frac{(f_\eps' F_\eps - F_\eps' f_\eps)^2}{(f_\eps^2 - F_\eps^2)^2} - \frac{2}{r^2} (1 - q_\eps)(3-q_\eps) \quad \textrm{ in } (0, 1].
	\label{Eq:S}
\end{align}

\begin{lemma}
Let $n \geq 2$, $\eps > 0$ and $\phi_\eps$ and $S_\eps$ be given by \eqref{Eq:phi} and \eqref{Eq:S}. Then
$$ -\nabla \cdot [(f_\eps^2 - F_\eps^2)\nabla \phi_\eps]-\frac{n-1}{r^2}  (f_\eps^2 - F_\eps^2)\phi_\eps= (f_\eps^2 - F_\eps^2)S_\eps \phi_\eps  \quad \textrm{ in } B \setminus \{0\},$$
that is, 
\begin{equation}
- \frac{[r^{n-1}(f_\eps^2 - F_\eps^2) \phi_\eps']'}{r^{n-1}(f_\eps^2 - F_\eps^2)}  - \frac{n-1}{r^2} \phi_\eps
= 
S_\eps \phi_\eps \quad  \textrm{ in } (0, 1) .
	\label{Eq:phi''}
\end{equation}
\end{lemma}

\begin{proof}
Let 
\begin{align}
b_\eps
	&= \frac{f_\eps f_\eps' - F_\eps F_\eps'}{f_\eps^2 - F_\eps^2} \quad \textrm{ in } (0, 1),\nonumber\\
\eta_\eps 
	&= \frac{\phi_\eps'}{\phi_\eps}
	= \frac{2}{r}q_\eps-\frac{2}{r}   -  b_\eps \quad \textrm{ in } (0, 1).
	\label{Eq:Sc-1}
\end{align}
We compute
\begin{align}
\frac{[r^{n-1}(f_\eps^2 - F_\eps^2) \phi_\eps']'}{r^{n-1}(f_\eps^2 - F_\eps^2)\phi_\eps} 
	&= \frac{[r^{n-1}(f_\eps^2 - F_\eps^2) \phi_\eps \eta_\eps]'}{r^{n-1}(f_\eps^2 - F_\eps^2) \phi_\eps}\nonumber\\
	&= \frac{[r^{n-3}(f_\eps^2 - F_\eps^2)^{1/2} F_\eps^2 \eta_\eps]'}{r^{n-3}(f_\eps^2 - F_\eps^2)^{1/2} F_\eps^2}
		= \eta_\eps' + \Big(\frac{2}{r} q_\eps + \frac{n-3}{r}  +  b_\eps\Big) \eta_\eps.
	\label{Eq:Sc-2}
\end{align}

We proceed to compute $\eta_\eps'$. First, since
\[
\frac{F_\eps''}{F_\eps} = \Big(\frac{F_\eps'}{F_\eps}\Big)' + \Big(\frac{F_\eps'}{F_\eps}\Big)^2 = \frac{q_\eps'}{r} + \frac{q_\eps^2 - q_\eps}{r^2},
\]
the ODE \eqref{Eq:Feps} for $F_\eps$ can be rewritten as
\begin{equation}
q_\eps' = - \frac{1}{r} (q_\eps - 1)(q_\eps + n - 1)  -\frac{r}{\eps^2}(1 - F_\eps^2).
	\label{Eq:q'}
\end{equation}
Next,
\[
b_\eps'
	= \frac{f_\eps f_\eps'' - F_\eps F_\eps''}{f_\eps^2 - F_\eps^2} + \frac{(f_\eps')^2 - (F_\eps')^2}{f_\eps^2 - F_\eps^2}- 2b_\eps^2.
\]
Using \eqref{Eq:feps} and \eqref{Eq:Feps}, we get
\[
b_\eps'
	= \frac{n-1}{r^2} - \frac{n-1}{r} b_\eps - \frac{1}{\eps^2}(1 - f_\eps^2 - F_\eps^2) + \frac{(f_\eps')^2 - (F_\eps')^2}{f_\eps^2 - F_\eps^2}- 2b_\eps^2.
\]
Therefore,  by \eqref{Eq:Sc-1},
\begin{align}
\eta_\eps' 
	&= - \frac{1}{r^2} (2q_\eps^2 + 2(n-1)q_\eps - (n+1))  
	+ \frac{n-1}{r} b_\eps  + 2b_\eps^2\nonumber\\
	&\qquad - \frac{1}{\eps^2}(1 + f_\eps^2 - F_\eps^2) 
	 - \frac{(f_\eps')^2 - (F_\eps')^2}{f_\eps^2 - F_\eps^2}.
	 \label{Eq:Sc-3}
\end{align}
Inserting \eqref{Eq:Sc-1} and \eqref{Eq:Sc-3} into \eqref{Eq:Sc-2} gives
\begin{align*}
\frac{[r^{n-1}(f_\eps^2 - F_\eps^2) \phi_\eps']'}{r^{n-1}(f_\eps^2 - F_\eps^2)\phi_\eps} 
	&=    \frac{1}{r^2} ( 2q_\eps^2 - 8q_\eps - (n-7))  
	- \frac{1}{\eps^2}(1 + f_\eps^2 - F_\eps^2) 
	 \\
	&\qquad 
	+  b_\eps^2  - \frac{(f_\eps')^2 - (F_\eps')^2}{f_\eps^2 - F_\eps^2}.
\end{align*}
Using the identity
\[
b_\eps^2  - \frac{(f_\eps')^2 - (F_\eps')^2}{f_\eps^2 - F_\eps^2} = \Big( \frac{f_\eps f_\eps' - F_\eps F_\eps'}{f_\eps^2 - F_\eps^2}\Big)^2  - \frac{(f_\eps')^2 - (F_\eps')^2}{f_\eps^2 - F_\eps^2} =  \frac{(f_\eps' F_\eps  - F_\eps' f_\eps)^2}{(f_\eps^2 - F_\eps^2)^2},
\]
we arrive at \eqref{Eq:phi''} after a simple rearrangement of terms.
\end{proof}

\begin{lemma}\label{Lem:CalS2}
Let $n \geq 2$ and $\eps > 0$. Then \eqref{Eq:tmcFS} holds.
\end{lemma}

\begin{proof} By \eqref{eq:ftilde}, $\widetilde{\mcF}_\eps$ is continuous in $H_0^1(B,\RR^n)$
because $F_\eps/f_\eps$ is continuous and positive in $[0,1]$. Consider first the case $w \in C_c^\infty(B\setminus \{0\},\RR^n)$.
Recalling that $w = f_\eps \tilde w$, we use the Hardy decomposition $\tilde w= \phi_\eps \hat w$ for the operator $-\nabla \cdot [(f_\eps^2 - F_\eps^2)\nabla]$ (see e.g. \cite[Lemma A.1]{INSZ3}) and \eqref{Eq:phi''} to find  
\begin{align*}
\int_B (f_\eps^2 - F_\eps^2) |\nabla \tilde w|^2\,dx
	&=  \int_B \Big\{(f_\eps^2 - F_\eps^2) \phi_\eps^2 |\nabla \hat w|^2  +|\hat w|^2  \phi_\eps \cdot \big[-\nabla \cdot \big((f_\eps^2 - F_\eps^2)  \nabla \phi_\eps\big)\big]\Big\}\,dx \\
&	= \int_B (f_\eps^2 - F_\eps^2) \Big[\phi_\eps^2 |\nabla \hat w|^2 + \Big(\frac{n-1}{r^2} + S_\eps\Big)\phi_\eps^2 |\hat w|^2  \Big]\,dx.
\end{align*}
Inserting this into \eqref{Eq:Cal4} we obtain \eqref{Eq:tmcFS} for $w\in C_c^\infty(B\setminus \{0\},\RR^n)$. 

Consider now the general case $w \in H_0^1(B,\RR^n)$. By density, we may select $\{w_j\} \subset C_c^\infty(B\setminus \{0\},\RR^n)$ such that $w_j \rightarrow w$ in $H_0^1(B,\RR^n)$. We have by the continuity of $\widetilde{\mcF}_\eps$ in $H_0^1(B,\RR^n)$ that
\begin{equation}
\widetilde{\mcF}_\eps(w) = \lim_{j \rightarrow \infty} \widetilde{\mcF}_\eps(w_j) = \lim_{j \rightarrow \infty}  \int_{B} (f_\eps^2 - F_\eps^2) \Big[ \phi_\eps^2\Big|\nabla \Big(\frac{w_j }{f_\eps \phi_\eps}\Big)\Big|^2 + \frac{S_\eps |w_j|^2}{f_\eps^2}\Big]  \,dx.
	\label{Eq:Fa0}
\end{equation}
Noting that $S_\eps = O_\eps(1)$ as $r \rightarrow 0$ (see \eqref{Eq:FAs}, \eqref{Eq:qFAs} and \eqref{Eq:feaAs}), we have that
\begin{align}
\lim_{j\rightarrow \infty}  \int_{B}\frac{ (f_\eps^2 - F_\eps^2) S_\eps |w_j|^2}{f_\eps^2}  \,dx 
	&=  \int_{B}\frac{ (f_\eps^2 - F_\eps^2) S_\eps |w|^2}{f_\eps^2}  \,dx,\label{Eq:Fa1}\\
\lim_{j,k\rightarrow \infty}  \int_{B}\frac{ (f_\eps^2 - F_\eps^2) S_\eps |w_j - w_k|^2}{f_\eps^2}  \,dx 
	&= 0.\label{Eq:Fa2}
\end{align}
By \eqref{Eq:Fa2} and the continuity of $\widetilde{\mcF}_\eps$, we have that
\begin{multline*}
\lim_{j,k\rightarrow \infty}  \int_{B} (f_\eps^2 - F_\eps^2) \phi_\eps^2 \Big|\nabla \Big(\frac{w_j - w_k}{f_\eps \phi_\eps}\Big)\Big|^2  \,dx\\
	 = \lim_{j,k\rightarrow \infty} \Big[\widetilde{\mcF}_{\eps}(w_j - w_k) -  \int_{B}\frac{ (f_\eps^2 - F_\eps^2) S_\eps |w_j - w_k|^2}{f_\eps^2}  \,dx \Big]
	= 0.
\end{multline*}
Hence the sequence $\Big\{(f_\eps^2 - F_\eps^2)^{1/2} \phi_\eps  \nabla \big(\frac{w_j}{f_\eps \phi_\eps}\big)\Big\}$ is Cauchy in $L^2(B, \RR^{n\times n})$ and hence converges therein to $(f_\eps^2 - F_\eps^2)^{1/2} \phi_\eps  \nabla \big(\frac{w}{f_\eps \phi_\eps}\big)$. In particular
\begin{equation}
\lim_{j \rightarrow \infty}  \int_{B} (f_\eps^2 - F_\eps^2) \phi_\eps^2 \Big|\nabla \Big(\frac{w_j  }{f_\eps \phi_\eps}\Big)\Big|^2  \,dx
	= \int_{B} (f_\eps^2 - F_\eps^2) \phi_\eps^2 \Big|\nabla \Big(\frac{w   }{f_\eps \phi_\eps}\Big)\Big|^2  \,dx.
	\label{Eq:Fa3}
\end{equation}
Using \eqref{Eq:Fa1} and \eqref{Eq:Fa3} in \eqref{Eq:Fa0}, we obtain \eqref{Eq:tmcFS}.
\end{proof}

In view of \eqref{Eq:Cal2}, \eqref{Eq:Cal3} and \eqref{Eq:tmcFS}, to prove the minimality of the vortex solution, it suffices to show that $S_\eps > 0$ in $(0,1)$. To this end, we interpolate $f_\eps$ and $F_\eps$ by the solutions $f_{\eps,a}$ of \eqref{Eq:fepsa} for $a \in [ F_\eps(1),1]$. By \cite{Hervex2, ODE_INSZ}, $f_{\eps,a}$ exists uniquely, $f_{\eps,a}' > 0$ and $f_{\eps,a}$ depends smoothly in $(\eps,a)$ and is increasing in $a$. Moreover, $f_\eps = f_{\eps,1}$ and $F_\eps = f_{\eps,F_\eps(1)}$ in $(0,1]$. 

For $a > F_\eps(1)$, define
\begin{equation}
S_{\eps,a} = \frac{1}{\eps^2}(1 + f_{\eps,a}^2 - F_\eps^2) - \frac{(f_{\eps,a}' F_\eps - F_\eps' f_{\eps,a})^2}{(f_{\eps,a}^2 - F_\eps^2)^2} - \frac{2}{r^2} (1 - q_\eps)(3-q_\eps) \quad \textrm{ in } (0,1).
	\label{Eq:Sepsa}
\end{equation}
Note that $S_\eps = S_{\eps,1}$. 

\begin{proposition}\label{Prop:Key}
Let $n \geq 2$ and $\eps > 0$. Then $\partial_a S_{\eps,a} > 0$ in $(0,1]$ for all $a >F_\eps(1)$ and $\lim_{a \rightarrow F_\eps(1)} S_{\eps,a} > 0$ in $(0,1]$. Consequently $S_{\eps,a} > 0$ in $(0,1]$ for all $a \in (F_\eps(1),1]$.
\end{proposition}

\begin{proof}
The result follows from Propositions \ref{Prop:pBnd} and  \ref{Prop:p<}.
\end{proof}

We are now ready to prove our main result.

\begin{proof}[Proof of Theorem \ref{Thm:Main}]
By Proposition \ref{Prop:Key}, $S_\eps = S_{\eps, 1} > 0$ in $(0,1]$. The result then follows from \eqref{Eq:Cal2}, \eqref{Eq:Cal3} and \eqref{Eq:tmcFS}.
\end{proof}

\section{Refined estimates for the entire vortex solution}\label{Sec3}

In this section, we gather properties of $F_\eps$ and associated quantities, culminating in a proof of the bound $\lim_{a \searrow F_\eps(1)} S_{\eps,a} > 0$. These estimates will also be used in proving the monotonicity in $a$ (and hence the positivity) of $S_{\eps,a}$ in Section \ref{Sec4}.

\subsection{Estimates for $F_\eps$ and $q_\eps$}

As is apparent from the expression of $S_{\eps,a}$ defined in \eqref{Eq:Sepsa}, in order to prove its positivity for $a > F_\eps(1)$, one needs to control $\frac{r^2}{\eps^2}$ from below in terms of $q_\eps$ defined in \eqref{Eq:q}. We establish:
\begin{proposition}\label{Prop:S31}
Let $n \geq 2$ and $\eps > 0$. Then $0 < q_\eps < 1$, $q_\eps' < 0$ and
\[
2q_\eps(1-q_\eps) < - rq_\eps' < \begin{cases}
	\frac{2q_\eps(1-q_\eps)(7 + 4q_\eps)}{4 + 7 q_\eps} &\text{ if } n = 2,\\
	3q_\eps(1-q_\eps) & \text{ if } n \ge 3
\end{cases} \text{ in } (0,\infty).
\]
Moreover,
\[
\frac{r^2}{\eps^2} > \frac{n+2}{2} (1 - q_\eps)(3 - q_\eps) \text{ in } (0,\infty).
\]
\end{proposition}

The proof of Proposition \ref{Prop:S31} is akin to some arguments in \cite[Section 4]{ODE_INSZ}, but new ideas are also needed.
We start with:

\begin{lemma}\label{Lem:qProp}
Let $n \geq 2$ and $\eps > 0$. Then $q_\eps(0) = 1$, $q_\eps'(0) = q_\eps(\infty) = q_\eps'(\infty) = 0$, $0 < q_\eps < 1$ and $q_\eps' < 0$ in $(0,\infty)$.
\end{lemma}

\begin{proof} A routine argument as in \cite{Hervex2, ODE_INSZ} shows that there exists a constant $d_0 > 0$ such that
\begin{equation}
\begin{cases}
F_\eps(r) = d_0 \eps^{-1} r - \frac{d_0}{2(n+2)} \eps^{-3}r^3 + \frac{d_0(2(n+2)d_0^2 + 1)}{8(n+2)(n+4)}\eps^{-5} r^5 + O(\eps^{-7} r^7) \text{ as } r \rightarrow 0,\\
F_\eps'(r) = d_0 \eps^{-1}  - \frac{3d_0}{2( n+2 )} \eps^{-3}r^2 +   \frac{5d_0(2(n+2)d_0^2 + 1)}{8(n+2)(n+4)}\eps^{-5} r^4 + O(\eps^{-7}r^6) \text{ as } r \rightarrow 0,\\
F_\eps(r) = 1 - \frac{n-1}{2} \eps^2 r^{-2} + \frac{3(n-1)(n-5)}{8} \eps^{4} r^{-4} + O(\eps^6 r^{-6}) \text{ as } r \rightarrow \infty, \\ 
F_\eps'(r) =   (n-1) \eps^2 r^{-3} - \frac{3(n-1)(n-5)}{2} \eps^{4} r^{-5} + O(\eps^6 r^{-7}) \text{ as } r \rightarrow \infty.
\end{cases}
\label{Eq:FAs}
\end{equation}
It follows that
\begin{equation}
\begin{cases}
q_\eps(r) = 1 - \frac{1}{n+2} \eps^{-2}r^2 + \frac{(n+2)^2 d_0^2 - 1}{(n+2)^2(n+4)} \eps^{-4} r^4 + O(\eps^{-6}r^6) \text{ as } r \rightarrow 0,\\
q_\eps'(r) = - \frac{2}{n+2} \eps^{-2}r  + \frac{4(n+2)^2 d_0^2 - 4}{(n+2)^2(n+4)} \eps^{-4} r^3 + O(\eps^{-6}r^5) \text{ as } r \rightarrow 0,\\
q_\eps(r) = (n-1) \eps^2 r^{-2} + O(\eps^{4} r^{-4}) \text{ as } r \rightarrow \infty, \\ 
q_\eps'(r) = -2(n-1) \eps^2 r^{-3} + O(\eps^{4} r^{-5}) \text{ as } r \rightarrow \infty.
\end{cases}
\label{Eq:qFAs}
\end{equation}
Hence $q_\eps(0) = 1$, $q_\eps'(0) = q_\eps(\infty) = q_\eps'(\infty) = 0$.

From the monotonicity of $F_\eps$, it follows by definition that $q_\eps > 0$. To see that $q_\eps < 1$, consider the function $\check F_\eps = F_\eps/r$ which satisfies
\[
\check F_\eps'' + \frac{n+1}{r} \check F_\eps' = - \frac{1}{\eps^2}(1 - F_\eps^2)\check F_\eps < 0.
\]
In particular, $r^{n+1}\check F_\eps'$ is a decreasing function. Since its limit as $r \rightarrow 0$ is zero, it follows that $\check F_\eps' < 0$ in $(0,\infty)$. This implies $\frac{F_\eps'}{F_\eps} - \frac{1}{r} = \frac{\check F_\eps'}{\check F_\eps} < 0$ and $q_\eps = \frac{rF_\eps'}{F_\eps} < 1$ in $(0,\infty)$.

It remains to show $q_\eps' < 0$ in $(0,\infty)$. Differentiating \eqref{Eq:q'} and using $r(F_\eps^2)' = 2 F_\eps^2 q_\eps$ we find
\begin{equation}
q_\eps'' = - \frac{1}{r} (2 q_\eps + n - 2)q_\eps' + \frac{1}{r^2} (q_\eps - 1)(q_\eps + n - 1) + \frac{1}{\eps^2} (2q_\eps + 1)F_\eps^2 - \frac{1}{\eps^2}.
	\label{Eq:q''}
\end{equation}
In particular,
\[
- r^2 q_\eps'' - r (2 q_\eps + n - 1)q_\eps' + \frac{2r^2F_\eps^2 }{\eps^2} q_\eps= \frac{2r^2}{\eps^2}(1 - F_\eps^2) > 0.
\]
Equation \eqref{Eq:q''} can be written without reference to $F_\eps$ as follows. First, by \eqref{Eq:q'},
\begin{equation}
\frac{F_\eps^2}{\eps^2} = \frac{q_\eps'}{r} + \frac{1}{r^2} (q_\eps - 1)(q_\eps + n - 1) + \frac{1}{\eps^2}.
	\label{Eq:q'X}
\end{equation}
Using \eqref{Eq:q'X} in \eqref{Eq:q''}, we get
\begin{equation}
q_\eps'' = - \frac{1}{r} (n - 3)q_\eps' +  \frac{2}{r^2} (q_\eps + 1)(q_\eps - 1)(q_\eps + n - 1) + \frac{2q_\eps}{\eps^2}.
	\label{Eq:q''X}
\end{equation}
Differentiating \eqref{Eq:q''X}, we get
\begin{align}
q_\eps''' &= - \frac{1}{r} (n - 3)q_\eps'' + \frac{1}{r^2} (6q_\eps^2 + 4(n-1)q_\eps + n-5) q_\eps' + \frac{2q_\eps'}{\eps^2}  \nonumber\\
	&\qquad -  \frac{4}{r^3} (q_\eps + 1)(q_\eps - 1)(q_\eps + n - 1)  .
	\label{Eq:q'''}
\end{align}
By \eqref{Eq:q'},
\[
- \frac{4}{r^3} (q_\eps+ 1) (q_\eps - 1)(q_\eps + n - 1) = \frac{4}{r^2} (q_\eps + 1) q_\eps'  + \frac{4}{\eps^2 r}(q_\eps + 1)(1 - F_\eps^2).
\]
Using this for the last term in \eqref{Eq:q'''} we get
\[
-q_\eps''' - \frac{1}{r} (n - 3)q_\eps'' + \Big[\frac{1}{r^2} (6q_\eps^2 + 4n q_\eps + n-1) + \frac{2}{\eps^2} \Big]q_\eps'
	= -\frac{4}{\eps^2 r}(q_\eps + 1)(1 - F_\eps^2) > 0.
\]
In particular, $q_\eps'$ cannot have a non-negative local maximum in $(0,\infty)$. As $q_\eps'(0) = q_\eps'(\infty) = 0$, we thus have $q_\eps' < 0$ in $(0,\infty)$. 
\end{proof}

\begin{figure}[h]
\begin{center}
\includegraphics[width=.8\textwidth]{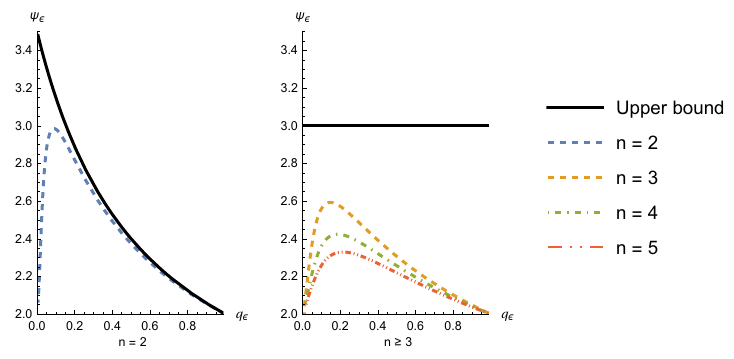}
\begin{minipage}{.8\textwidth}
\caption{An illustration of the bounds for $\psi_\eps$. Here the plotted graphs are those of $\psi_\eps \circ q_\eps^{-1}$. The black curves indicate the respective upper bounds in different dimensions in Lemma \ref{Lem:psiUB}.}
\label{Fig1}
\end{minipage}
\end{center}
\end{figure}

It is convenient to introduce
\[
\psi_\eps = \frac{rq_\eps'}{q_\eps( q_\eps - 1)} > 0
\]
so that the estimates in Proposition \ref{Prop:S31} become
\[
2 < \psi_\eps <  \begin{cases}
\frac{2(7+4q_\eps)}{4+ 7 q_\eps} & \text{ if } n = 2,\\
3 & \text{ if } n \geq 3
\end{cases} \quad 
\text{ in }(0,\infty) .
\]
See Figure \ref{Fig1} for an illustration.

\begin{lemma}\label{Lem:psiLB}
Let $n \geq 2$ and $\eps > 0$. Then $\psi_\eps(0) = \psi_\eps(\infty) = 2$ and $\psi_\eps > 2$ in $(0,\infty)$.
\end{lemma}

\begin{proof} By \eqref{Eq:qFAs}, $\psi_\eps(0) = \psi_\eps(\infty) = 2$. Differentiating $\psi_\eps$, we find
\[
\psi_\eps' = \frac{q_\eps''}{q_\eps'} \psi_\eps  - \frac{1}{r}(2q_\eps - 1) \psi_\eps^2+ \frac{1}{r} \psi_\eps.
\]
Using \eqref{Eq:q''X}, we get
\begin{equation}
\psi_\eps' = - \frac{1}{r} (n - 4)  \psi_\eps - \frac{1}{r}(2q_\eps - 1) \psi_\eps^2 + \frac{2}{r q_\eps} (q_\eps + 1)(q_\eps + n - 1)  - \frac{2r}{\eps^2(1 - q_\eps)} .
	\label{Eq:psi'}
\end{equation}
We rewrite this as
\[
\frac{1}{r} (1 - q_\eps) \psi_\eps' = \frac{1}{r^2} A(q_\eps,\psi_\eps) - \frac{2}{\eps^2}
\]
with
\[
A(x,y) = (1-x)\Big[- (n - 4)  y - (2x - 1) y^2 + \frac{2}{x} (x + 1)(x + n - 1)\Big].
\]
Differentiating and using $q_\eps' = \frac{1}{r} \psi_\eps q_\eps (q_\eps - 1)$, we get 
\begin{equation}
(1 - q_\eps) \psi_\eps''
	- \frac{1}{r} [(1 - q_\eps) (1 - \psi_\eps q_\eps ) + \partial_y A(q_\eps,\psi_\eps)]\psi_\eps' 
	= \frac{1}{r^2} \frac{1-q_\eps}{q_\eps} B(q_\eps, \psi_\eps)
	\label{Eq:psi''}
\end{equation}
with
\begin{align*}
B(x,y) 
	&= - x^2y \partial_x A(x,y) - \frac{2x}{1-x} A(x,y)\\
	&= - x^2 y^3(4x - 3) - xy^2[(n-8)x + 2] \\
		&\qquad + 2y[2x^3 + (n-1)x^2 + (n-4)x + (n-1)] \\
		&\qquad- 4(x+1)(x+n-1).
\end{align*}

By \eqref{Eq:psi''}, at an interior minimum point of $\psi_\eps$, it holds that
\[
B(q_\eps, \psi_\eps) \geq 0.
\]
Therefore, as $\psi_\eps(0) = \psi_\eps(\infty) = 2$, to show that $\psi_\eps > 2$, we only need to show that
\begin{equation}
B(x,y) < 0 \text{ for } x \in (0,1), y \in (0,2].
	\label{Eq:Bneg}
\end{equation}
Indeed, note that $B$ is non-increasing with respect to dimension $n$ for $x \in [0,1], y \in (0,2]$:
\[
\frac{\partial B}{\partial n} = - (2-y)(2 + 2x - x^2y) \leq 0.
\]
Thus, we may take $n = 2$ when proving \eqref{Eq:Bneg}. Let $B_2$ denote $B$ when $n = 2$, i.e.
\[
B_2(x,y) 
	= x^2 y^3(3 - 4x) - 2xy^2(1-3x) 
		+ 2y(2x^3 + x^2 - 2x + 1)
			- 4(x+1)^2.
\]
Note that
\[
B_2(x,2) 
	= -24x(x-1)^2 < 0 \text{ for } x \in (0,1).
\]
Therefore, \eqref{Eq:Bneg} is proved once we show that $B_2$ is increasing in $y$ for $x \in [0,1], y \in (0,2]$. We compute
\[
\partial_y B_2(x,y) = 3x^2 y^2(3 - 4x) - 4xy(1-3x) + 2\big(2x^3 + (x-1)^2\big).
\]
It is clear that $\partial_y B_2(x,y) > 0$ for $x \in (1/3,3/4)$, $y \in [0,2]$ (as all three summands are positive). When $0 \leq x \leq 1/3$, we have
\[
3 - 4x \geq \frac{5}{3} \text{ and } 2(2x^3 + x^2 - 2x + 1) \geq 2(x-1)^2 \geq \frac{8}{9} > \frac{4}{5}, 
\]
and so
\[
\partial_y B_2(x,y) > 5x^2y^2 - 4xy + \frac{4}{5} = 5(xy - \frac{2}{5})^2 \geq 0.
\]
When $x \geq 3/4$, $\partial_y B_2(x,\cdot)$ is concave in $y$, so its minimum in $[0,2]$ is attained at the endpoints. For $x\in [0,1]$,
\begin{align*}
\partial_y B_2(x,0) &= 4x^3 + 2(x-1)^2 > 0,\\
\partial_y B_2(x,2) &= -44x^3 + 62 x^2 - 12 x + 2  =44(x^2-x^3)+  2(3x - 1)^2 > 0 ,
\end{align*}
yielding $\partial_y B_2(x,y) > 0$ for $x \in [3/4,1]$, $y \in [0,2]$. 
In any case, we have shown that $B_2$ is increasing in $y$ for $x \in [0,1], y \in (0,2]$. As explained before, this proves \eqref{Eq:Bneg} and the inequality $\psi_\eps > 2$ in $(0,\infty)$.
\end{proof}

\begin{lemma}\label{Lem:psiUB}
Let $n \geq 2$ and $\eps > 0$. Then 
\[
\psi_\eps <  \begin{cases}
\frac{2(7+4q_\eps)}{4+ 7 q_\eps} & \text{ if } n = 2,\\
3 & \text{ if } n \geq 3
\end{cases}
\quad \text{ in }(0,\infty) .
\]
\end{lemma}

The upper bound $\psi_\eps < 3$ in dimensions $n \geq 3$ can be sharpened to $\psi_\eps < 3-q_\eps$, but we decide not to include it here as we do not need it. In dimension $n = 2$, it follows from our estimate above that $\psi_\eps < 7/2$. However, this global bound would not be strong enough for our later use in the region $r \approx 0$, i.e. $q_\eps \approx 1$, where we will use that $\psi_\eps \approx 2$.

\begin{proof}
Define
\[
g(x) = \begin{cases}
\frac{2x(1 - x) (7 + 4 x)}{4 + 7 x} & \text{ if } n = 2,\\
3x(1-x) & \text{ if } n \geq 3.
\end{cases}
\]
We need to show that
\[
\chi_\eps= rq_\eps' + g(q_\eps) > 0 \text{ in } (0,\infty).
\]
Suppose by contradiction that this does not hold. As $\chi_\eps(0) = \chi_\eps(\infty) = 0$, this implies that $\chi_\eps$ has a non-positive interior minimum.

We compute
\[
\chi_\eps'
	= r q_\eps'' + (1 + g'(q_\eps)) q_\eps'.
\]
Using \eqref{Eq:q''X}, we get
\begin{align}
\chi_\eps'
	&=   (g'(q_\eps) - n + 4) q_\eps'
	+ \frac{2}{r} (q_\eps + 1)(q_\eps - 1)(q_\eps + n - 1) + \frac{2rq_\eps}{\eps^2} 	\nonumber\\
	&=   \frac{1}{r}(g'(q_\eps) - n + 4) \chi_\eps\nonumber\\
		&\qquad
		+ \frac{1}{r}[2(q_\eps + 1)(q_\eps - 1)(q_\eps + n - 1) - (g'(q_\eps) - n + 4)g(q_\eps) ]
	  + \frac{2rq_\eps}{\eps^2}.
\label{Eq:chi'}
\end{align}
We rewrite this as
\[
\frac{1}{rq_\eps} \chi_\eps'
	= \frac{1}{r^2} M(q_\eps,  \chi_\eps) + \frac{2}{\eps^2} 
\]
where
\begin{align*}
M(x,z) 
	&= c_1(x) z + c_0(x),\\
c_1(x)
	&= \frac{1}{x}(g'(x) - n + 4),\\
c_0(x)
	&= \frac{1}{x}\Big[2(x + 1)(x - 1)(x + n - 1) - (g'(x) - n + 4)g(x)\Big].
\end{align*}
Differentiating and using $rq_\eps' = \chi_\eps - g(q_\eps) $, we get
\begin{equation}
\frac{1}{q_\eps} \chi_\eps'' = \frac{1}{rq_\eps^2}[q_\eps + rq_\eps'+ q_\eps^2 \partial_z M (q_\eps, \chi_\eps)]\chi_\eps' + \frac{1}{r^2} P(q_\eps,\chi_\eps)
	\label{Eq:chi''}
\end{equation}
where
\begin{align*}
P(x,z) 
	&= \partial_x M(x,z)(z - g(x))  - 2M(x,z)\\
	&= [c_1'(x) z + c_0'(x)][z - g(x)] - 2[c_1(x) z + c_0(x)]\\
	&= c_1'(x) z^2 + [c_0'(x) - c_1'(x) g(x) - 2 c_1(x)]z - c_0'(x) g(x) - 2 c_0(x).
\end{align*}
In particular, 
\begin{equation}
P(q_\eps,\chi_\eps) \geq 0 \text{ at an interior minimum of } \chi_\eps.
	\label{Eq:Psign}
\end{equation}

We continue the proof with the case $3 \leq n \leq 7$, followed by the case $n \geq 8$, and then $n = 2$.

\medskip
\noindent\underline{Case 1:} $3 \leq n \leq 7$.

In this case, $g(x) = 3x (1 - x)$ and 
\begin{align}
P(x,z) 
	&= \frac{1}{x^2} \Big[(n - 7) z^2  + a_1(x) z + x(x-1) a_0(x)\Big],\label{Eq:P3-1}\\
a_1(x)
	&= -32x^3 + 2(n + 14) x^2   - (n-7)x + 2(n-1),\label{Eq:P3-2}\\
a_0(x)
	&=- 96x^3
			- (3n - 143)x^2  + 2(n -21)x
			+ 2 (n-1).\label{Eq:P3-3}
\end{align}

In light of \eqref{Eq:Psign} and the statement that $\chi_\eps$ has a non-positive interior minimum, we only need to show that
\begin{equation}
P(x,z) < 0 \text{ for all } x \in (0,1), z \leq 0.
	\label{Eq:Pcond}
\end{equation}
To this end, it is enough to show that
\begin{equation}
a_1(x) > 0 \text{ and } a_0(x) > 0 \text{ for } x \in [0,1].
	\label{Eq:a01cond}
\end{equation}
Indeed, for $a_1$, we have
\begin{align*}
a_1(x)
	&= 2(n-1)(1 - x^3) + (7-n)(x - x^3) \\
		&\qquad + 2(n+14)(x^2 - x^3) + (3n + 1)x^3 > 0 \text{ for } x \in [0,1].
\end{align*}
For $a_0$, note that $a_0$ is increasing with respect to $n$ for $x \in (0,1)$:
\[
\frac{\partial a_0}{\partial n} =  2(1-x)(2x + 1) > 0.
\]
Therefore, it suffices to show that
\[
a_{0,3}(x):= a_0(x)\Big|_{n = 3} = - 2(48 x^3 - 67x^2 + 18 x - 2) > 0 \text{ in } [0,1].
\]
A simple minimization gives 
\[
\min_{[0,1]} a_{0,3} = a_{0,3}\Big(\frac{1}{144} (67 - 1897^{1/2})\Big) = \frac{1}{15552}(102475 - 1897^{3/2}) \approx 1.276 > 0.
\]
This proves \eqref{Eq:a01cond} and completes the proof when $3 \leq n \leq 7$.

\medskip
\noindent\underline{Case 2:} $n\geq 8$. 

Formulas \eqref{Eq:P3-1}-\eqref{Eq:P3-3} remain valid. The term $(n - 7)$ in front of $z^2$ in the expression of $P$ now has an unfavorable sign for the argument in Case 1. To circumvent this, observe from \eqref{Eq:chi'} that $M(q_\eps, \chi_\eps) < 0$ at a critical point of $\chi_\eps$. Together with \eqref{Eq:Psign}, we have that 
\begin{equation}
\tilde P(q_\eps,\chi_\eps) \geq P(q_\eps,\chi_\eps) \geq 0 \text{ at a non-positive interior minimum of } \chi_\eps,
	\label{Eq:PsignX}
\end{equation}
where
\[
\tilde P(x,z) = P(x,z) + \frac{1}{x} (M(x,z) + 6z) z = \frac{1}{x }\big[ \tilde a_1 (x) z +  (x-1)a_0(x)\big]
\]
with
\begin{align*}
\tilde a_1(x) 
	&= -48x^2 + (n + 65)x + 2(n-8)\\
	& = 2(n-8)(1-x^2) + (n+65)(x-x^2) + (3n + 1)x^2 > 0.
\end{align*}
As $a_0 \geq a_{0,3} > 0$ in $[0,1]$ (from the argument in Case $1$), it is readily seen that $\tilde P(x,z) < 0$ for all $x \in (0,1), z \leq 0$, and therefore \eqref{Eq:PsignX} cannot hold. This completes the proof when $n \geq 8$.

\medskip
\noindent\underline{Case 3:} $n = 2$. 

In this case, $g(x) = \frac{2x (1 - x)(7 + 4 x)}{4 + 7x}$ and 
\begin{align*}
P(x,z) 
	&= \frac{1}{x^2} \Big[-\frac{88}{(4+7x)^3} b_2(x) z^2  + \frac{2}{(4+7x)^4} b_1(x) z + \frac{12x(x-1)^3}{(4+7x)^5}  b_0(x)\Big], \\
b_2(x)
	&=  4 + 21 x + 12 x^2 + 7 x^3 ,  \\
b_1(x)
	&=256 + 2848 x + 23088 x^2 + 10144 x^3 + 10609 x^4\\
		&\qquad + 9120 x^5 + 3969 x^6 - 1470 x^7, \\
b_0(x)
	&=256 - 1920 x + 5800 x^2 + 8920 x^3 + 1265 x^4 - 4375 x^5 - 1960 x^6. 
\end{align*}
To reach a contradiction to \eqref{Eq:Psign}, it suffices to show that  
\begin{equation}
b_2(x) > 0, b_1(x) > 0 \text{ and } b_0(x) > 0 \text{ for } x \in [0,1].
	\label{Eq:b012cond}
\end{equation}
It is clear that $b_2 > 0$ in $[0,1]$. Also, since $x^7 \leq x^6$, we have that $b_1 > 0$ in $[0,1]$. For $b_0$, using $x^3 \geq x^5$ and $x^3 \geq x^6$, we have
\[
b_0(x) \geq 256 - 1920 x + 5800 x^2 > 0.
\]
This proves \eqref{Eq:b012cond}, which gives a contradiction to \eqref{Eq:Psign}. The conclusion in the case $n = 2$ follows.
\end{proof}

We give a few consequences of the established bounds for $\psi_\eps$.

\begin{corollary}\label{Cor:F<>q}
Let $n \geq 2$ and $\eps > 0$. Then $1 - F_\eps^2  < q_\eps$ in $(0,\infty)$.
\end{corollary}

\begin{proof}
Let $\beta = \frac{1 - F_\eps^2}{q_\eps}$. By \eqref{Eq:FAs} and \eqref{Eq:qFAs}, we have $\beta(0) = \beta(\infty) = 1$. Since $F_\eps' = \frac{1}{r} F_\eps q_\eps$ and $q_\eps' = - \frac{1}{r} \psi_\eps q_\eps (1 - q_\eps)$,
\[
\beta' 
	= \frac{-2F_\eps^2 q_\eps  + (1-F_\eps^2)  \psi_\eps   (1 - q_\eps)}{rq_\eps }
	=  \frac{1}{r}\big[-2 + \beta(2q_\eps + \psi_\eps(1 - q_\eps))\big].
\]
As $q_\eps < 1$ (by Lemma \ref{Lem:qProp}) and $\psi_\eps > 2$ (by Lemma \ref{Lem:psiLB}), this implies
\[
\beta' > \frac{1}{r}(-2 + 2\beta).
\]
In particular, $\beta < 1$ at its critical points. As $\beta(0) = \beta(\infty) = 1$, this implies $\beta < 1$ everywhere and so $1 - F_\eps^2 < q_\eps$ in $(0,\infty)$.
\end{proof}

\begin{corollary}\label{Cor:G2}
Let $n \geq 2$ and $\eps > 0$. Then $\frac{(1 - q_\eps)(3q_\eps + n - 1)}{q_\eps} < \frac{r^2}{\eps^2}$ in $(0,\infty)$.
\end{corollary}

\begin{proof}
By \eqref{Eq:q'} and the fact that $\psi_\eps > 2$ (by Lemma \ref{Lem:psiLB}),
\[
\frac{r^2}{\eps^2}(1 - F_\eps^2) = - r q_\eps' + (1 - q_\eps)(q_\eps + n - 1) > (1 - q_\eps)(3q_\eps + n - 1). 
\]
The conclusion follows from the estimate $1 - F_\eps^2 < q_\eps$ in Corollary \ref{Cor:F<>q}.
\end{proof}

\begin{corollary}\label{Cor:rq'>stuff}
Let $n \geq 2$ and $\eps > 0$. Then $rq_\eps' + (1 - q_\eps)(q_\eps + n - 1) > 0$ in $(0,\infty)$.
\end{corollary}

\begin{proof}
The conclusion follows from Lemma \ref{Lem:psiUB}, the identity $rq_\eps' = - \psi_\eps q_\eps (1 - q_\eps)$, and the elementary inequality
\[
x + n - 1 > \begin{cases}
	\frac{2x(7 + 4x)}{4 + 7x} & \text{ if } n = 2,\\
	3x  &\text{ if } n \geq 3
\end{cases}
	\quad  \text{ for } 0 < x < 1.
\]
\end{proof}

\begin{lemma}\label{Lem:G3/2}
Let $n \geq 2$ and $\eps > 0$. Then $(1 - q_\eps)(3-q_\eps) < \frac{2}{n+2} \frac{r^2}{\eps^2}$ in $(0,\infty)$.
\end{lemma}

Note that, when $n \geq 4$, the result follows from Corollary \ref{Cor:G2} and the fact that
\[
\max_{x\in [0,1]} \frac{x(3 - x)}{3x + n - 1} = \begin{cases}
\frac{11-2\sqrt{10}}{9} &\text{ if } n = 2,\\
	\frac{13-2\sqrt{22}}{9}  &\text{ if } n = 3,\\
	\frac{2}{n+2}& \text{ if } n \geq 4.
\end{cases}
\]
We give a different argument which works for all dimensions $n \geq 2$. We will use the following Gronwall-type inequality, whose simple proof we skip.

\begin{lemma}[Gronwall's inequality]\label{Lem:Gronwall}
Let $\theta \in C^1((0,r_0))$ and $a \in C^0([0,r_0])$  for some $r_0 > 0$, and assume that
\[
\theta'(r) + \frac{1}{r} a(r) \theta(r) > 0 \text{ in } (0,r_0).
\]
If 
\[
\lim_{r \rightarrow 0} \theta(r) r^{a(0)} \geq 0 \quad \text{ and }\quad a(r) - a(0) = O(r^\alpha) \text{ for some } \alpha > 0 \text{ near } r = 0,
\]
then $\theta(r) > 0$ in $(0,r_0)$.
\end{lemma}

\begin{proof}[Proof of Lemma \ref{Lem:G3/2}]
Let
\[
\xi = \frac{r^2}{\eps^2} - \frac{n+2}{2} (1 - q_\eps)(3 - q_\eps).
\]
We compute using \eqref{Eq:q'}
\begin{align*}
\xi'
	&= \frac{2r}{\eps^2} + \frac{n+2}{r} (2- q_\eps)(1-q_\eps)(q_\eps + n - 1) - \frac{(n+2)r}{\eps^2}(2-q_\eps)(1 - F_\eps^2)\\
	&= \frac{2}{r}\Big[1 - \frac{n+2}{2}(2-q_\eps)(1 - F_\eps^2)\Big] \xi \\
		&\quad 
		+ \frac{n+2}{r}(1 - q_\eps) \Big\{(3 - q_\eps) \Big[1 - \frac{n+2}{2}(2-q_\eps)(1 - F_\eps^2)\Big] + (2- q_\eps)(q_\eps + n - 1) \Big\}.
\end{align*}
Using the bound $1 - F_\eps^2 < q_\eps$ (see Corollary \ref{Cor:F<>q}), we obtain
\begin{align*}
\xi'
	&> \frac{2}{r}\Big[1 - \frac{n+2}{2}(2-q_\eps)(1 - F_\eps^2)\Big] \xi \\
		&\quad 
		+ \frac{n+2}{r}(1 - q_\eps) \Big\{(3 - q_\eps) \Big[1 - \frac{n+2}{2}(2-q_\eps)q_\eps\Big] + (2- q_\eps)(q_\eps + n - 1) \Big\}\\
	&= \frac{2}{r}\Big[1 - \frac{n+2}{2}(2-q_\eps)(1 - F_\eps^2)\Big] \xi \\
		&\quad 
		+ \frac{n+2}{2r}(1 - q_\eps)^2 \big[2(2n+1) - 2(2n+3) q_\eps + (n+2) q_\eps^2\big].
\end{align*}
As
\[
2(2n+1) - 2(2n+3) q_\eps + (n+2) q_\eps^2 > 0 \text{ for } 0 < q_\eps < 1,
\]
it follows that 
\[
\xi' > \frac{2}{r}\Big[1 - \frac{n+2}{2}(2-q_\eps)(1 - F_\eps^2)\Big] \xi.
\]
By Gronwall's inequality (see Lemma \ref{Lem:Gronwall}), $\xi > 0$ in $(0,\infty)$. Equivalently, $(1 - q_\eps)(3-q_\eps) < \frac{2}{n+2} \frac{r^2}{\eps^2}$ in $(0,\infty)$.
\end{proof}

\begin{proof}[Proof of Proposition \ref{Prop:S31}]
The result follows from Lemmas \ref{Lem:qProp}, \ref{Lem:psiLB}, \ref{Lem:psiUB} and \ref{Lem:G3/2}.
\end{proof}

\subsection{Estimates for the linearized ODE}

In this subsection, we prove $\lim_{a \searrow F_\eps(1)} S_{\eps,a} > 0$. Introducing
\[
t_\eps = \Big[\frac{r^2}{\eps^2} - 2(1 - q_\eps)(3-q_\eps)\Big]^{1/2} \quad \textrm{ in } (0, \infty),
\]
which is well-defined thanks to Lemma \ref{Lem:G3/2}, we have
\[
\lim_{a \searrow F_\eps(1)} S_{\eps,a} = \frac{1}{r^2} t_\eps^2   - \lim_{a \searrow F_\eps(1)} \frac{(f_{\eps,a}' F_\eps - F_\eps' f_{\eps,a})^2}{(f_{\eps,a}^2 - F_\eps^2)^2}  \quad \textrm{ in } (0, 1].
\]
Let $h_\eps$ denote the linearization of  \eqref{Eq:Feps}:
\[
\begin{cases}
h_{\eps}'' + \frac{n - 1}{r} h_{\eps}' - \frac{n-1}{r^2} h_{\eps} = -\frac{1}{\eps^2} (1 - 3F_{\eps}^2)h_{\eps} \text{ in } (0,\infty),\\
h_\eps(0) = 0, h_\eps (1) = 1.
\end{cases}
\]
It is known that such $h_\eps$ exists uniquely, and in fact, $h_\eps = \partial_a f_{\eps,a}\Big|_{a = F_{\eps}(1)}$ in $(0,1]$. In particular, $h_\eps > 0$ and $h_\eps'(0) > 0$. Defining
\[
p_\eps = \frac{r h_\eps'}{h_\eps} \quad \textrm{ in } (0, \infty), 
\]
we find by L'H\^opital's rule that
\[
\lim_{a \searrow F_\eps(1)} \frac{(f_{\eps,a}' F_\eps - F_\eps' f_{\eps,a})^2}{(f_{\eps,a}^2 - F_\eps^2)^2} = \frac{1}{4r^2} (p_\eps - q_\eps)^2.
\]
Therefore
\begin{equation}
\lim_{a \searrow F_\eps(1)} S_{\eps,a} = \frac{1}{r^2} t_\eps^2   - \frac{1}{4r^2} (p_\eps - q_\eps)^2 \quad \textrm{ in } (0, 1].
	\label{Eq:Slowend}
\end{equation}

\begin{proposition}\label{Prop:pBnd}
Let $n \geq 2$ and $\eps > 0$. Then $q_\eps < p_\eps < q_\eps + 2t_\eps$ in $(0,\infty)$. In particular, $\lim_{a \searrow F_\eps(1)} S_{\eps,a} > 0$ in $(0,1]$.
\end{proposition}

\begin{figure}[ht]
\begin{center}
\includegraphics[width=.5\textwidth]{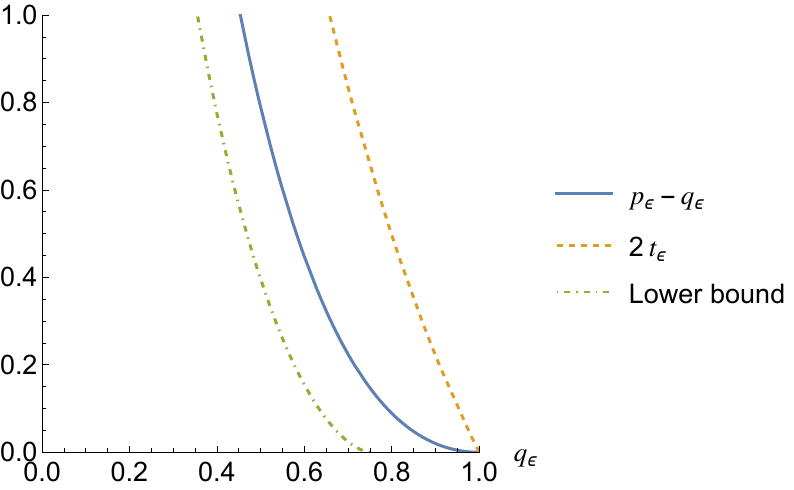}
\begin{minipage}{.6\textwidth}
\caption{An illustration of the upper bound in Proposition \ref{Prop:pBnd} and the lower bound in Lemma \ref{Lem:pF>Sn2} for $n = 2$. All functions are plotted as functions of $q_\eps$.}
\label{Fig2}
\end{minipage}
\end{center}
\end{figure}

See Figure \ref{Fig2} for an illustration.

By \eqref{Eq:FAs} and the ODE for $h_\eps$, it follows 
\[
\begin{cases}
h_\eps(r) = h_\eps'(0)   r - \frac{h_\eps'(0)}{2(n+2)} \eps^{-2} r^3 + O_\eps(r^5) \text{ as } r \rightarrow 0,\\
h_\eps'(r) = h_\eps'(0)  - \frac{3h_\eps'(0)}{2(n+2)} \eps^{-2} r^2 + O_\eps(r^4) \text{ as } r \rightarrow 0,
\end{cases}
\]
which implies
\begin{equation}
p_\eps(r) = 1 - \frac{1}{ n+2} \eps^{-2} r^2 + O_\eps(r^4) \text{ as } r \rightarrow 0.
	\label{Eq:pFAs}
\end{equation}

We start with the lower bound for $p_\eps$.

\begin{lemma}\label{Lem:pF>}
Let $n \geq 2$ and $\eps > 0$. Then $p_\eps > q_\eps$ in $(0,\infty)$.
\end{lemma}

\begin{proof}
We have
\begin{equation}
p_{\eps}' = - \frac{1}{r} (p_{\eps} - 1)(p_{\eps} + n - 1) - \frac{r}{\eps^2}(1 - 3 F_{\eps}^2).
	\label{Eq:pF'}
\end{equation}
Combining with \eqref{Eq:q'}, we get
\[
(p_{\eps} - q_{\eps})' + \frac{1}{r}(p_{\eps} + q_{\eps} + n - 2) (p_{\eps} - q_{\eps}) 
	= \frac{2r}{\eps^2} F_{\eps}^2  > 0.
\]
In view of \eqref{Eq:qFAs} and \eqref{Eq:pFAs}, we have by Gronwall's inequality (see Lemma \ref{Lem:Gronwall}) that $p_{\eps} > q_{\eps}$ in $(0,\infty)$.
\end{proof}

In Section \ref{Sec4}, we will need the following strengthened lower estimate for $p_\eps$ in dimension $n = 2$ away from $r = 0$. (It is not needed for the proof of Proposition \ref{Prop:pBnd}.) See Figure \ref{Fig2} for an illustration.

\begin{lemma}\label{Lem:pF>Sn2}
Let $n = 2$ and $\eps > 0$. Let $r_1 > 0$ be such that $q_\eps(r_1) = \frac{3}{4}$. Then
\[
p_\eps > q_\eps + \frac{7r}{4\eps}  \big(\frac{3}{4} - q_\eps\big)^{3/2}  \text{ in } [r_1,\infty).
\]
\end{lemma}

\begin{proof}
We claim that the functions
\[
\mu(r) := \frac{7r}{4\eps}\big(\frac{3}{4} - q_\eps(r)\big)^{3/2} \text{ for } r \ge r_1
\]
satisfy
\begin{equation}
\mu' + \frac{1}{r} \mu^2 + \frac{2q_\eps}{r} \mu - \frac{2r}{\eps^2} (1-q_\eps) < 0 \text{ in } (r_1, \infty).
\label{Eq:ClaimMu}
\end{equation}
Postponing the proof of this claim for the moment, let us show that it implies the conclusion. Indeed, by \eqref{Eq:q'}, \eqref{Eq:pF'}, the fact that $1 - q_\eps < F_\eps^2$ (see Corollary \ref{Cor:F<>q}), and the claim, we have
\begin{equation}
(p_\eps - q_\eps - \mu)' + \frac{1}{r}(p_\eps + q_\eps + \mu )(p_\eps - q_\eps - \mu) > 0 \text{ in } [r_1,\infty).
	\label{Eq:pqmu}
\end{equation}
As $p_\eps - q_\eps > 0$ (by Lemma \ref{Lem:pF>}) and $\mu(r_1) = 0$, we have by Gronwall's inequality that 
\[
p_\eps - q_\eps - \mu > 0 \text{ in } [r_1,\infty).
\]

It remains to prove the claimed differential inequality \eqref{Eq:ClaimMu}. 

\medskip
\noindent\underline{Proof of \eqref{Eq:ClaimMu}.}
Using the fact that $rq_\eps' = - \psi_\eps q_\eps(1 - q_\eps)$, we compute
\begin{align*}
&\mu' + \frac{1}{r} \mu^2 + \frac{2q_\eps }{r} \mu - \frac{2r}{\eps^2} (1-q_\eps)\\
	&\qquad =  \frac{7 (2q_\eps +  1)}{4\eps}\big(\frac{3}{4} - q_\eps \big)^{3/2}  + \frac{21}{8\eps}\big(\frac{3}{4} - q_\eps \big)^{1/2}  \psi_\eps q_\eps(1-q_\eps) 
		 \\
		 &\qquad\qquad - \frac{2 r}{\eps^2}\big[(1-q_\eps) - \frac{49}{32} \big(\frac{3}{4} - q_\eps \big)^3 		
		    \big].
\end{align*}
Observe that, as $0 < q_\eps \leq \frac{3}{4} $, we have
\[
(1-q_\eps) - \frac{49}{32} \big(\frac{3}{4} - q_\eps \big)^3 \geq (1-q_\eps) - \frac{49}{32} \big(\frac{3}{4} - q_\eps \big)^2.
\]
The right hand side is a concave quadratic in $q_\eps$ and its values when $q_\eps = 0$ and $q_\eps = \frac{3}{4}$ are positive. Hence
\[
(1-q_\eps) - \frac{49}{32} \big(\frac{3}{4} - q_\eps \big)^3 > 0.
\]
 Therefore, by the estimates $\frac{r^2}{\eps^2} > \frac{(1-q_\eps)(3q_\eps + 1)}{q_\eps}$ (see Corollary \ref{Cor:G2}) and $\psi_\eps < \frac{1}{2}(7 - 3q_\eps)$ (see Lemma \ref{Lem:psiUB}), we have that 
\begin{align*}
&\mu' + \frac{1}{r} \mu^2 + \frac{2q_\eps }{r} \mu - \frac{2r}{\eps^2} (1 - q_\eps)\\
	&\qquad \leq  \frac{7 (2q_\eps +  1)}{4\eps}\big(\frac{3}{4} - q_\eps \big)^{3/2}  + \frac{21}{16\eps}\big(\frac{3}{4} - q_\eps \big)^{1/2}  q_\eps(1-q_\eps) (7 - 3q_\eps)
		 \\
		 &\qquad\qquad - \frac{2}{\eps}\frac{(1-q_\eps)^{1/2} (3q_\eps + 1)^{1/2}}{q_\eps^{1/2}} \big[(1-q_\eps) - \frac{49}{32} \big(\frac{3}{4} - q_\eps \big)^3 		
		    \big].
\end{align*}
Therefore, to prove \eqref{Eq:ClaimMu}, we only need to check the algebraic inequality
\begin{align}
&   \big(\frac{3}{4} - x \big)^{3/2} + \frac{3 }{4 }\big(\frac{3}{4} - x \big)^{1/2}   \frac{x(1-x)(7-3x)}{2x+1}   \nonumber\\
		&\quad -   \frac{8  (1-x)^{1/2} (3x + 1)^{1/2}}{7x^{1/2}(2x+1)}\big[(1-x) - \frac{49}{32}\big(\frac{3}{4} - x \big)^3 \big] < 0 \text{ for } 0 \leq x \leq \frac{3}{4}.
		\label{Eq:AE-1}
\end{align}
As the proof of this inequality is tedious, although elementary, we move it to the appendix. See Lemma \ref{Lem:Alg1}.
\end{proof}

\begin{lemma}\label{Lem:pF<}
Let $n \geq 2$ and $\eps > 0$. Then $p_\eps < q_\eps + 2t_\eps$ in $(0,\infty)$.
\end{lemma}

\begin{proof}
We first claim that
\begin{equation}
E :=  2t_\eps^2 + \frac{r^2}{\eps^2 t_\eps}[1 - t_\eps - (2 - q_\eps - t_\eps)(1-F_\eps^2)] > 0 \text{ in } (0,\infty).
	\label{Eq:pqt-1}
\end{equation}
We split the argument into two cases.

\medskip
\noindent\underline{Case 1:} $2 - q_\eps - t_\eps \leq 0$. 

In this case, using $F_\eps < 1$, followed by $t_\eps < \frac{r}{\eps}$,
\[
E >  2t_\eps^2  + \frac{r^2}{\eps^2 t_\eps}(1 - t_\eps) >   2 t_\eps^2 + t_\eps - \frac{r^2}{\eps^2} = t_\eps^2 + t_\eps - 2 (1 - q_\eps)(3-q_\eps) = t_\eps^2 + t_\eps + 2 - 2(2 - q_\eps)^2.
\]
As $t_\eps \geq 2 - q_\eps$ (by assumption) and $0 < q_\eps < 1$, it follows that
\[
E \geq  3q_\eps -   q_\eps^2 > 0.
\]
\medskip
\noindent\underline{Case 2:} $2 - q_\eps - t_\eps > 0$. 

In this case, using $1 - F_\eps^2 < q_\eps$ (see Corollary \ref{Cor:F<>q}), 
\begin{align*}
E &>  2t_\eps^2 + \frac{r^2}{\eps^2 t_\eps}\big[1 - t_\eps - (2 - q_\eps - t_\eps)q_\eps\big]\\
	&= 2t_\eps^2 + \frac{r^2}{\eps^2 t_\eps}(1- q_\eps)(1 - q_\eps - t_\eps).
\end{align*}
If $t_\eps \leq 1 - q_\eps$, the above inequality implies $E > 0$. Assume henceforth that $t_\eps > 1 - q_\eps$. Substituting in $\frac{r^2}{\eps^2} = t_\eps^2 + 2(1 - q_\eps)(3-q_\eps)$, we find
\begin{align*}
E 
	&> \frac{1}{t_\eps} \big[(1 + q_\eps) t_\eps^3 + (1 - q_\eps)^2 t_\eps^2 - 2 (1 - q_\eps)^2(3-q_\eps) t_\eps + 2(1-q_\eps)^3(3-q_\eps)\big].
\end{align*}
Letting $\hat t_\eps = \frac{t_\eps}{1 - q_\eps} > 1$ and using $\hat t_\eps^3 > \hat t_\eps^2$, we obtain
\begin{align*}
E 
	&> \frac{(1-q_\eps)^2}{\hat t_\eps} \big[(1 + q_\eps) \hat t_\eps^3 + (1 - q_\eps)  \hat t_\eps^2 - 2  (3-q_\eps) \hat t_\eps + 2(3-q_\eps)\big]\\
	&> \frac{2(1-q_\eps)^2}{\hat t_\eps} \big[ \hat t_\eps^2 -  (3-q_\eps) \hat t_\eps +  (3-q_\eps)\big]\\
	&>  2(1-q_\eps)^2 \big[2 \sqrt{3-q_\eps} -  (3-q_\eps)\big] > 0.
\end{align*}
The claim is proved.

We compute using \eqref{Eq:q'}
\begin{align}
t_\eps' 
	&= \frac{1}{t_\eps} \Big( \frac{r}{\eps^2} + 2(2 - q_\eps )q_\eps'\Big)\nonumber\\
	&= \frac{1}{ r t_\eps} \Big[(2 - q_\eps )r q_\eps' + (2 - q_\eps )(1 - q_\eps)(q_\eps + n - 1) + \frac{r^2}{\eps^2} - \frac{r^2}{\eps^2}(2 - q_\eps)(1 - F_\eps^2) \Big].
	\label{Eq:t'-1}
\end{align}
Observe that \eqref{Eq:pqt-1} is equivalent to
\[
\frac{1}{r t_\eps}\Big[\frac{r^2}{\eps^2} - \frac{r^2}{\eps^2}(2 - q_\eps)(1 - F_\eps^2)\Big] > -\frac{2}{r} t_\eps^2 + \frac{r}{\eps^2} F_\eps^2.
\]
Using this estimate together with Corollary \ref{Cor:rq'>stuff} in \eqref{Eq:t'-1} gives
\[
t_\eps' >  - \frac{2}{r} t_\eps^2  +  \frac{r}{\eps^2} F_\eps^2 .
\]
Combining with \eqref{Eq:q'}, we find
\[
(q_\eps + 2t_\eps)'
	> - \frac{1}{r}(q_\eps + 2 t_\eps - 1)(q_\eps + 2t_\eps + n - 1) - \frac{r}{\eps^2}(1 - 3F_\eps^2).
\]
Recalling \eqref{Eq:pF'}, we obtain
\[
(p_\eps - q_\eps - 2t_\eps)' + \frac{1}{r}(p_\eps + q_\eps + 2t_\eps + n - 2)(p_\eps - q_\eps - 2t_\eps) < 0.
\]
By Gronwall's inequality (see Lemma \ref{Lem:Gronwall}), we deduce that $p_\eps - q_\eps - 2t_\eps < 0$.
\end{proof}

\begin{proof}[Proof of Proposition \ref{Prop:pBnd}]
The result follows from \eqref{Eq:Slowend} together with Lemmas \ref{Lem:pF>} and \ref{Lem:pF<}.
\end{proof}

\section{Monotonicity and positivity of $S_{\eps,a}$}
\label{Sec4}

In this section, we prove that $\partial_a S_{\eps,a} > 0$. Throughout the section, we assume that 
\[
a > F_\eps(1).
\] 
Parallel to $q_\eps$, $h_\eps$ and $p_\eps$, for $f_{\eps,a}$ defined at \eqref{Eq:fepsa}, we introduce
\[
q_{\eps,a} = \frac{r f_{\eps,a}'}{f_{\eps,a}}, \quad h_{\eps,a} = \partial_a f_{\eps,a} > 0, 
 \quad \text{ and } \quad  p_{\eps,a} = \frac{rh_{\eps,a}'}{h_{\eps,a}} \quad \textrm{ in } (0,1].
\]
It is routine to check that, with $d_{\eps,a} = f_{\eps,a}'(0) = O_\eps(1)$,
\begin{equation}
\begin{cases}
f_{\eps,a}(r) = d_{\eps,a} r - \frac{d_{\eps,a}}{2(n+2)} \eps^{-2}r^3 + \frac{d_{\eps,a}(2(n+2)d_{\eps,a}^2\eps^2  + 1)}{8(n+2)(n+4)}\eps^{-4} r^5 + O_\eps( r^7) \text{ as } r \rightarrow 0,\\
f_{\eps,a}'(r) = d_{\eps,a}  - \frac{3d_{\eps,a}}{2(n+2)} \eps^{-2}r^2+   \frac{5d_{\eps,a}(2(n+2)d_{\eps,a}^2\eps^2  + 1)}{8(n+2)(n+4)}\eps^{-4} r^4 + O_\eps( r^6) \text{ as } r \rightarrow 0,
\end{cases}
\label{Eq:feaAs}
\end{equation}
\begin{equation}
\begin{cases}
q_{\eps,a}(r) = 1 - \frac{1}{n+2} \eps^{-2}r^2 + \frac{ (n+2)^2d_{\eps,a}^2 \eps^2- 1}{(n+2)^2(n+4)} \eps^{-4}r^4 + O_\eps(r^6) \text{ as } r \rightarrow 0,\\
q_{\eps,a}'(r) = - \frac{2}{n+2} \eps^{-2}r  + \frac{4 (n+2)^2d_{\eps,a}^2 \eps^2- 4}{(n+2)^2(n+4)} \eps^{-4}r^3 + O_\eps(r^5) \text{ as } r \rightarrow 0,
\end{cases}
\label{Eq:qAs}
\end{equation}
and
\begin{equation}
p_{\eps,a}(r) = 1 - \frac{1}{n+2} \eps^{-2} r^2 + O_\eps(r^4) \text{ as } r \rightarrow 0.
	\label{Eq:pAs}
\end{equation}

We compute
\begin{align*}
\frac{1}{2}\partial_a S_{\eps,a} 
	&= \frac{1}{\eps^2} f_{\eps,a} h_{\eps,a} - \frac{(f_{\eps,a}' F_\eps - F_\eps' f_{\eps,a})(h_{\eps,a}' F_\eps - F_\eps' h_{\eps,a})}{(f_{\eps,a}^2 - F_\eps^2)^2} + \frac{2(f_{\eps,a}' F_\eps - F_\eps' f_{\eps,a})^2 f_{\eps,a} h_{\eps,a}}{(f_{\eps,a}^2 - F_\eps^2)^3}\\
	&= \frac{1}{\eps^2} f_{\eps,a} h_{\eps,a}\Big[1
		- \frac{\eps^2}{r^2} F_\eps^2\frac{(q_{\eps,a} - q_\eps)(p_{\eps,a} - q_\eps)}{(f_{\eps,a}^2 - F_\eps^2)^2}
		 + \frac{2\eps^2}{r^2} F_\eps^2 f_{\eps,a}^2 \frac{(q_{\eps,a} - q_\eps)^2}{(f_{\eps,a}^2 - F_\eps^2)^3}\Big].
\end{align*}
In particular,
\[
\partial_a S_{\eps,a}  > 0 \Leftrightarrow 
p_{\eps,a} <  q_\eps 
		 +  \frac{2 f_{\eps,a}^2(q_{\eps,a} - q_\eps)}{f_{\eps,a}^2 - F_\eps^2}
		 +
\frac{r^2}{\eps^2 F_\eps^2}\frac{(f_{\eps,a}^2 - F_\eps^2)^2}{q_{\eps,a} - q_\eps}.
\]
For convenience, we introduce the following scalar radial functions in $(0,1]$:
\begin{align*}
x_{\eps,a} 
	&= \frac{q_{\eps,a} - q_\eps}{f_{\eps,a}^2 - F_\eps^2},\\
y_{\eps,a}
	& = f_{\eps,a}^2 x_{\eps,a},\\
z_{\eps,a} 
	&= \frac{r^2}{\eps^2} \frac{f_{\eps,a}^2 - F_\eps^2}{F_\eps^2 x_{\eps,a}},\\
P_{\eps,a}
	&= q_\eps + 2y_{\eps,a} + z_{\eps,a}.
\end{align*}
Then $\partial_a S_{\eps,a}  > 0$ is equivalent to $p_{\eps,a} < P_{\eps,a}$ in $(0,1]$.

\begin{proposition}\label{Prop:p<}
Let $n \geq 2$, $\eps > 0$ and $a>F_\eps(1)$. Then $x_{\eps,a}, y_{\eps,a}, z_{\eps,a}, P_{\eps,a}$ are continuous in $[0,1]$, positive away from $r = 0$, and $q_\eps + 2y_{\eps,a} < p_{\eps,a} < P_{\eps,a}$ in $(0,1]$. In particular, $\partial_a S_{\eps,a} > 0$ in $(0,1]$.
\end{proposition}

Now, observe that $x_{\eps,a}$ and $y_{\eps,a}$ are continuous in $[0,1]$ (recalling $f_{\eps,a} > F_\eps$, $q_{\eps,a}(0) = q_\eps(0) = 1$ and $q_{\eps,a}'(0) = q_\eps'(0) = 0$ by \eqref{Eq:qFAs} and \eqref{Eq:qAs}). Note also that, by \eqref{Eq:qFAs} and \eqref{Eq:qAs}, $\lim_{r \rightarrow 0} \frac{x_{\eps,a}}{r^2} > 0$, which implies that $z_{\eps,a}$ is continuous in a neighborhood of $r = 0$. The fact that $z_{\eps,a}$ is continuous in $[0,1]$ is then given by the following lemma.
\begin{lemma}\label{Lem:qqp}
Let $n \geq 2$ and $\eps > 0$. Then $p_{\eps,a} > q_{\eps,a} > q_\eps$ in $(0,1]$.
\end{lemma}

\begin{proof}
Recall \eqref{Eq:q'}
\begin{align*}
q_{\eps}' = - \frac{1}{r}(q_\eps - 1)(q_\eps + n - 1) - \frac{r}{\eps^2}(1- F_\eps^2).
\end{align*}
Similarly,
\begin{align*}
q_{\eps,a}' = - \frac{1}{r}(q_{\eps,a} - 1)(q_{\eps,a} + n - 1) - \frac{r}{\eps^2}(1- f_{\eps,a}^2).
\end{align*}
Thus
\begin{equation}
(q_{\eps,a} - q_\eps)' + \frac{1}{r}(q_{\eps,a} + q_\eps + n - 2) (q_{\eps,a} - q_\eps) 
	= \frac{r}{\eps^2} (f_{\eps,a}^2 - F_\eps^2) > 0.
	\label{Eq:q-q'}
\end{equation}
By Gronwall's inequality (Lemma \ref{Lem:Gronwall}) together with \eqref{Eq:qFAs} and \eqref{Eq:qAs}, this implies $q_{\eps,a} - q_\eps > 0$ in $(0,1]$. 

Next, we have
\[
h_{\eps,a}'' + \frac{n - 1}{r} h_{\eps,a} - \frac{n-1}{r^2} h_{\eps,a} = -\frac{1}{\eps^2} (1 - 3f_{\eps,a}^2)h_{\eps,a}.
\]
Therefore
\begin{equation}
p_{\eps,a}' = - \frac{1}{r} (p_{\eps,a} - 1)(p_{\eps,a} + n - 1) - \frac{r}{\eps^2}(1 - 3 f_{\eps,a}^2).
	\label{Eq:p'}
\end{equation}
As before, this leads to
\[
(p_{\eps,a} - q_{\eps,a})' + \frac{1}{r}(p_{\eps,a} + q_{\eps,a} + n - 2) (p_{\eps,a} - q_{\eps,a}) 
	= \frac{2r}{\eps^2} f_{\eps,a}^2  > 0.
\]
By Gronwall's inequality (Lemma \ref{Lem:Gronwall}) together with  \eqref{Eq:qAs} and \eqref{Eq:pAs}, this implies $p_{\eps,a} > q_{\eps,a}$ in $(0,1]$.
\end{proof}

In the next two lemmas, we will show $p_{\eps,a} > q_\eps + 2y_{\eps,a}$. This inequality stresses the importance of the term $z_{\eps,a}$ in the inequality $p_{\eps,a} < P_{\eps,a} = q_\eps + 2y_{\eps,a} + z_{\eps,a}$ that we are trying to prove. We will only use it in the case $n = 2$.

\begin{lemma}\label{Lem:xy<}
Let $n \geq 2$ and $\eps > 0$. Then $2x_{\eps,a} y_{\eps,a} < \frac{r^2}{\eps^2}$ in $(0,1]$.
\end{lemma}

\begin{proof}
Since $F_\eps' = \frac{1}{r} F_\eps q_\eps$ and $f_{\eps,a}' = \frac{1}{r} f_{\eps,a}q_{\eps,a}$, we have
\begin{equation}
\frac{(f_{\eps,a}^2 - F_\eps^2)'}{f_{\eps,a}^2 - F_\eps^2}
	= \frac{2}{r} \frac{f_{\eps,a}^2 q_{\eps,a} - F_\eps^2 q_\eps}{f_{\eps,a}^2 - F_\eps^2}\\
	= \frac{2}{r} q_{\eps} + \frac{2}{r} y_{\eps,a}.
	\label{Eq:prex'}
\end{equation}
Hence, by \eqref{Eq:q-q'},
\begin{align}
\frac{x_{\eps,a}'}{x_{\eps,a}}
	&= \frac{(q_{\eps,a} - q_\eps)'}{q_{\eps,a} - q_\eps} 
		-\frac{(f_{\eps,a}^2 - F_\eps^2)'}{f_{\eps,a}^2 - F_\eps^2} \nonumber\\
	&= - \frac{1}{r}(q_{\eps,a} + 3q_\eps + n - 2) 
		 - \frac{2}{r} y_{\eps,a}
		  + \frac{r}{\eps^2 x_{\eps,a}} ,
	\label{Eq:x'}
\end{align}
and
\begin{align}
\frac{y_{\eps,a}'}{y_{\eps,a}}
	&=  \frac{2}{r} q_{\eps,a} + \frac{x_{\eps,a}'}{x_{\eps,a}}\nonumber \\
	&= - \frac{1}{r}(-q_{\eps,a} + 3q_\eps + n - 2) 
		 - \frac{2}{r} y_{\eps,a}
		  + \frac{r}{\eps^2 x_{\eps,a}}.
	\label{Eq:y'}
\end{align}
It follows that
\begin{equation}
(x_{\eps,a} y_{\eps,a})' 
	= -\frac{2}{r}(3q_\eps + n - 2 + 2y_{\eps,a}) x_{\eps,a} y_{\eps,a} 
	 +   \frac{2r}{\eps^2}y_{\eps,a}.
	 \label{Eq:xy'}
\end{equation}

By \eqref{Eq:xy'},
\begin{align*}
\Big(\frac{r^2}{\eps^2} - 2x_{\eps,a} y_{\eps,a}\Big)' 
	&=  \frac{4}{r}(3q_\eps + n - 2 + 2 y_{\eps,a}) x_{\eps,a} y_{\eps,a} 
	 +   \frac{2r}{\eps^2}(1 - 2y_{\eps,a})\\
	 &= \frac{2}{r} (1 - 2y_{\eps,a})\Big(\frac{r^2}{\eps^2} - 2x_{\eps,a} y_{\eps,a}\Big)
	 	+ \frac{4}{r}(3q_\eps + n - 1) x_{\eps,a} y_{\eps,a} \\
	&> \frac{2}{r} (1 - 2y_{\eps,a})\Big(\frac{r^2}{\eps^2} - 2x_{\eps,a} y_{\eps,a}\Big).
\end{align*}
Using that $x_{\eps,a}, y_{\eps,a}$ are continuous in $[0,1]$ and $x_{\eps,a}(r) = O(r^2)$ and $y_{\eps,a}(r) = O(r^4)$ near $r = 0$, we may apply Gronwall's inequality (Lemma \ref{Lem:Gronwall}) to obtain $\frac{r^2}{\eps^2} > 2x_{\eps,a} y_{\eps,a} $ in $(0,1]$.
\end{proof}

\begin{lemma}\label{Lem:p>}
Let $n \geq 2$ and $\eps > 0$. Then $p_{\eps,a} > q_\eps + 2y_{\eps,a}$ in $(0,1]$.
\end{lemma}

\begin{proof} By \eqref{Eq:y'},
\begin{align}
(q_\eps + 2y_{\eps,a})'
	&= - \frac{1}{r}(q_\eps + 2y_{\eps,a} - 1)(q_\eps + 2y_{\eps,a} + n - 1)\nonumber\\
		&\qquad 
		- \frac{r}{\eps^2}(1 - 3f_{\eps,a}^2) -  \frac{r}{\eps^2}(f_{\eps,a}^2 - F_\eps^2) \Big(1 - \frac{2\eps^2}{r^2} x_{\eps,a}y_{\eps,a}\Big).
	\label{Eq:q2y'}
\end{align}
Combining with \eqref{Eq:p'} gives
\begin{align*}
&[p_{\eps,a} - (q_\eps + 2y_{\eps,a})]'
	+\frac{1}{r} [p_{\eps,a} + (q_\eps + 2y_{\eps,a}) + n-2] [p_{\eps,a} - (q_\eps + 2y_{\eps,a})]\nonumber\\
		&\qquad = 
		  \frac{r}{\eps^2}(f_{\eps,a}^2 - F_\eps^2) \Big(1 - \frac{2\eps^2}{r^2} x_{\eps,a}y_{\eps,a}\Big) > 0,
\end{align*}
where we have used Lemma \ref{Lem:xy<} in the last inequality. 
In view of \eqref{Eq:qFAs}, \eqref{Eq:pAs} and the fact that $y_{\eps,a} = O(r^4)$ near $r = 0$, we have by Gronwall's inequality (Lemma \ref{Lem:Gronwall}) that $p_{\eps,a} - (q_\eps + 2y_{\eps,a}) > 0$ in $(0,1]$.
\end{proof}

\begin{lemma}\label{Lem:y>}
Let $n = 2$ and $\eps > 0$. Let $r_1 > 0$ be such that $q_\eps(r_1) = \frac{3}{4}$. Then
\[
y_{\eps,a} > \frac{7}{8} \frac{(1-q_\eps)^{1/2} (3q_\eps + 1)^{1/2}}{q_\eps^{1/2}}  \big(\frac{3}{4} - q_\eps\big)^{3/2}  \text{ in } [r_1,\infty).
\]
\end{lemma}

\begin{proof}
We compute
\begin{align*}
\partial_a q_{\eps,a} 
	&= \frac{r h_{\eps,a}'}{f_{\eps,a}} - \frac{r f_{\eps,a}'}{f_{\eps,a}^2} h_{\eps,a} 
		= \frac{h_{\eps,a}}{f_{\eps,a}}(p_{\eps,a} - q_{\eps,a}),\\
\partial_a x_{\eps,a}
	&= \frac{1}{f_{\eps,a}^2 - F_\eps^2} \frac{h_{\eps,a}}{f_{\eps,a}}(p_{\eps,a} - q_{\eps,a}) - \frac{2(q_{\eps,a} - q_\eps)}{(f_{\eps,a}^2 - F_\eps^2)^2} f_{\eps,a} h_{\eps,a}\\	&= \frac{h_{\eps,a}}{f_{\eps,a} (f_{\eps,a}^2 - F_\eps^2)} (p_{\eps,a} - q_{\eps,a} - 2 y_{\eps,a}),\\
\partial_a y_{\eps,a}
	&= 2 f_{\eps,a} h_{\eps,a} x_{\eps,a} +  f_{\eps,a}^2 \partial_a x_{\eps,a}\\
	&= \frac{f_{\eps,a} h_{\eps,a}}{(f_{\eps,a}^2 - F_\eps^2)} \big[p_{\eps,a} - q_{\eps,a} -  2 y_{\eps,a} + 2(f_{\eps,a}^2 - F_\eps^2)x_{\eps,a}\big].
\end{align*}
Recalling that $(f_{\eps,a}^2 - F_\eps^2)x_{\eps,a} = q_{\eps,a} - q_\eps$, we thus have
\[
\partial_a y_{\eps,a}
	= \frac{f_{\eps,a} h_{\eps,a}}{(f_{\eps,a}^2 - F_\eps^2)} \big[(p_{\eps,a} - q_\eps -   2y_{\eps,a}) + (q_{\eps,a} - q_\eps) \big].
\]
In particular, as $h_{\eps,a} > 0$, we have by the inequalities $q_{\eps,a} > q_\eps$ (see Lemma \ref{Lem:qqp}) and $p_{\eps,a}> q_\eps + 2y_{\eps,a}$ (see Lemma \ref{Lem:p>}) that
\[
\partial_a y_{\eps,a}
	 > 0.
\]
On the other hand, since
\[
\lim_{a \rightarrow F_\eps(1)} x_{\eps,a} = \frac{\frac{h_{\eps,a}}{f_{\eps,a}}(p_{\eps,a} - q_{\eps,a})}{2 f_{\eps,a} h_{\eps,a}} \Big|_{a = F_{\eps}(1)} = \frac{p_\eps - q_\eps}{2F_\eps^2}
\]
we have
\[
\lim_{a \rightarrow F_\eps(1)} y_{\eps,a} = \frac{1}{2}(p_\eps - q_\eps).
\]
Therefore,
\[
y_{\eps,a} > \frac{1}{2} (p_\eps - q_\eps).
\]
The conclusion follows from the lower bound for $\frac{r}{\eps}$ in Corollary \ref{Cor:G2} and the lower bound for $p_\eps$ in Lemma \ref{Lem:pF>Sn2}. 
\end{proof}

\begin{lemma}\label{Lem:xy>}
Let $n \geq 2$ and $\eps > 0$. Then $\frac{r^2}{\eps^2} < 2(3q_\eps + n - 1 + y_{\eps,a})x_{\eps,a}$ in $(0,1]$.
\end{lemma}

\begin{proof}
Consider
\[
\zeta = \frac{r^2}{\eps^2} - 2(3q_\eps + n - 1 + y_{\eps,a})x_{\eps,a}.
\]
By \eqref{Eq:x'} and the fact that $q_{\eps,a} = q_\eps + y_{\eps,a} - F_\eps^2 x_{\eps,a}$,
\begin{align*}
[(3q_\eps + n - 1)x_{\eps,a}]'
	&= \frac{r}{\eps^2}(3q_\eps + n - 1)
		 \\
		&\qquad
		 - \frac{1}{r}(3q_\eps + n - 1)(4q_\eps + n - 2 + 3y_{\eps,a} )x_{\eps,a}
		 + 3q_\eps' x_{\eps,a}\\
		 &\qquad
		 + \frac{1}{r}(3q_\eps + n - 1) F_\eps^2 x_{\eps,a}^2.
\end{align*}
Together with \eqref{Eq:xy'}, this gives
\begin{align*}
\zeta'
	& = - \frac{2r}{\eps^2} (3q_\eps + n - 2 +  2 y_{\eps,a}) 
		+ \frac{2}{r}(3q_\eps + n - 1)(4q_\eps + n - 2)x_{\eps,a}
	\nonumber \\
		&\qquad + \frac{2 }{r}(15 q_\eps + 5n - 7 + 4 y_{\eps,a}) x_{\eps,a} y_{\eps,a}  
		 - 6q_\eps' x_{\eps,a}
		 - \frac{2}{r}(3q_\eps + n - 1) F_\eps^2 x_{\eps,a}^2
		 \nonumber\\
	& = - \frac{2}{r} (3q_\eps + n - 2 +  2 y_{\eps,a}) \zeta
		- \frac{2}{r} (3q_\eps + n - 1 )(2q_\eps + n - 2) x_{\eps,a}
	\nonumber \\
		&\qquad - \frac{2 }{r}(3 q_\eps + n-1) x_{\eps,a} y_{\eps,a}  
		 - 6q_\eps' x_{\eps,a}
		 - \frac{2}{r}(3q_\eps + n - 1) F_\eps^2 x_{\eps,a}^2.
\end{align*}

Recalling that $q_\eps' = -\frac{1}{r}\psi_\eps q_\eps(1-q_\eps)$, we may use the upper bound of $\psi_\eps$ in Lemma \ref{Lem:psiUB} to proceed.

\medskip
\noindent\underline{Case 1:} $n \geq 3$. 

By the bound $\psi_\eps< 3$ in Lemma \ref{Lem:psiUB} when $n \geq 3$,
\[
- q_\eps' =  \frac{1}{r} \psi_\eps q_\eps(1 - q_\eps) < \frac{3}{r} q_\eps(1 - q_\eps).
\]
Therefore
\begin{align*}
\zeta'
	&< - \frac{2}{r} (3q_\eps + n - 2 +  2  y_{\eps,a}) \zeta\\
		&\qquad - \frac{2}{r}(15 q_\eps^2 + (5n - 17)q_\eps + (n-1)(n-2))x_{\eps,a}.
\end{align*}
Noting that
\[
15 q_\eps^2 + (5n - 17)q_\eps + (n-1)(n-2) > 0 \text{ for } n \geq 3,
\]
we then have
\[
\zeta'
	< - \frac{2}{r} (3q_\eps + n - 2 +  2  y_{\eps,a}) \zeta \text{ in } (0,\infty).
\]
Since $\zeta = O(r^2)$, $q_\eps = 1 + O(r^2)$ and $ y_{\eps,a} = O(r^4)$ near $r = 0$, we have by Gronwall's inequality (Lemma \ref{Lem:Gronwall}) that $\zeta < 0$ in $(0,1]$, as desired.

\medskip
\noindent\underline{Case 2:} $n = 2$.

Arguing as in Case 1, but using the bound $\psi_\eps < \frac{2(7 + 4q_\eps)}{4 + 7 q_\eps} < \frac{1}{2}(7 - 3q_\eps)$ (see Lemma \ref{Lem:psiUB}) when $n = 2$, we get
\begin{align*}
\zeta'
	&< - \frac{2}{r} (3q_\eps   +  2  y_{\eps,a}) \zeta  - \frac{1}{r}\big[q_\eps(-9 q_\eps^2 + 42q_\eps  - 17) + 2(3q_\eps + 1)y_{\eps,a} \big]x_{\eps,a}.
\end{align*}
Note that, in $[0,1]$, the quadratic polynomial $-9t^2 + 42t - 17$ is positive if and only if $t > \frac{1}{3} (7 - 4 \sqrt{2}) \approx .44$. Let $r_0$ be such that $q_\eps(r_0) = \frac{1}{3} (7 - 4 \sqrt{2})$, then 
\begin{equation}
\zeta'
	< - \frac{2}{r} (3q_\eps   +  2  y_{\eps,a}) \zeta \text{ in } [0,r_0].
	\label{Eq:zetaP1}
\end{equation}
On the other hand, in $(r_0, \infty)$, we have $q_\eps < \frac{1}{3} (7 - 4 \sqrt{2}) < \frac{3}{4}$ and we may use the lower bound of $y_{\eps,a}$ from Lemma \ref{Lem:y>} to get
\begin{align*}
\zeta'
	&< - \frac{2}{r} (3q_\eps  +  2  y_{\eps,a}) \zeta\\
		&\qquad - \frac{1}{r}\Big[-q_\eps(9 q_\eps^2 - 42q_\eps  + 17) + \frac{7}{4}   \frac{(1-q_\eps)^{1/2} (3q_\eps + 1)^{3/2}}{q_\eps^{1/2}}  \big(\frac{3}{4} - q_\eps\big)^{3/2} \Big]x_{\eps,a}.
\end{align*}
By Lemma \ref{Lem:Alg2} in the appendix (noting that $\frac{1}{3} (7 - 4 \sqrt{2}) < \frac{1}{2}$), the second line of the above inequality is negative. Thus
\[
\zeta'
	< - \frac{2}{r} (3q_\eps   +  2  y_{\eps,a}) \zeta \text{ in } (r_0,\infty).
\]
Recalling also \eqref{Eq:zetaP1}, we may apply Gronwall's inequality (Lemma \ref{Lem:Gronwall}) to conclude as in Case 1.
\end{proof}

\begin{proof}[Proof of Proposition \ref{Prop:p<}] By Lemma \ref{Lem:p>}, $p_{\eps,a} > q_\eps + 2y_{\eps,a}$. It remains to show that $p_{\eps,a} < P_{\eps,a} = q_\eps + 2y_{\eps,a} + z_{\eps,a}$. By \eqref{Eq:prex'} and \eqref{Eq:x'},
\begin{align*}
z_{\eps,a}'
	&= z_{\eps,a}\Big(\frac{2}{r} + \frac{(f_{\eps,a}^2 - F_\eps^2)'}{f_{\eps,a}^2 - F_\eps^2} - \frac{2}{r} q_\eps - \frac{x_{\eps,a}'}{x_{\eps,a}} \Big)\\
	&=  \frac{1}{r} \Big(q_{\eps,a} + 3q_\eps + n 
		 + 4 y_{\eps,a}
		  - \frac{r^2}{\eps^2 x_{\eps,a}}  \Big)z_{\eps,a}.
\end{align*}
Combining with \eqref{Eq:q2y'} and noting that
\[
\frac{r}{\eps^2}(f_{\eps,a}^2 - F_\eps^2) \Big(1 - \frac{2\eps^2}{r^2} x_{\eps,a}y_{\eps,a}\Big)
	= \frac{1}{r} F_{\eps}^2 z_{\eps,a} x_{\eps,a} \Big(1 - \frac{2\eps^2}{r^2} x_{\eps,a}y_{\eps,a}\Big),
\]
we get
\begin{align*}
P_{\eps,a}'
	&= - \frac{1}{r}(P_{\eps,a} - 1)(P_{\eps,a} + n - 1)
		- \frac{r}{\eps^2}(1 - 3f_{\eps,a}^2)\\
		&\quad + \frac{1}{r} \Big(q_{\eps,a} + 5q_\eps + 2n - 2 
		 + 8 y_{\eps,a} + z_{\eps,a}
		 + \frac{2\eps^2}{r^2} F_{\eps}^2x_{\eps,a}^2y_{\eps,a}
		 - \frac{r^2}{\eps^2 x_{\eps,a}} - F_{\eps}^2 x_{\eps,a} \Big)z_{\eps,a}.
\end{align*}
In particular, since $q_{\eps,a} > q_\eps$ (see Lemma \ref{Lem:qqp}), $F_\eps^2 x_{\eps,a} < f_{\eps,a}^2 x_{\eps,a} = y_{\eps,a}$, and 
\[
\frac{r^2}{\eps^2 x_{\eps,a}} < 2(3q_\eps + n - 1 + y_{\eps,a}) \quad (\text{by Lemma \ref{Lem:xy>}}),
\]
we have
\begin{align*}
P_{\eps,a}'
	&> - \frac{1}{r}(P_{\eps,a} - 1)(P_{\eps,a} + n - 1)
		- \frac{r}{\eps^2}(1 - 3f_{\eps,a}^2).
\end{align*}
Recalling \eqref{Eq:p'}, we obtain
\[
(p_{\eps,a} - P_{\eps,a})'
	+ \frac{1}{r} (p_{\eps,a} + P_{\eps,a} + n - 2) (p_{\eps,a} - P_{\eps,a}) < 0.
\]
By Gronwall's inequality (Lemma \ref{Lem:Gronwall}) together with \eqref{Eq:qFAs}, \eqref{Eq:pAs}, and the fact that $y_{\eps,a} = O(r^4)$ and $z_{\eps,a} = O(1)$ near $r = 0$, we have that $p_{\eps,a} < P_{\eps,a}$ in $(0,1]$.
\end{proof}

\appendix
\section{Some algebraic inequalities of one variable}

In this appendix, we give the proof of the two algebraic inequalities \eqref{Eq:AE-1} and \eqref{Eq:AEX-1} that were used in the body of the paper.

\begin{lemma}\label{Lem:Alg1}
The inequality \eqref{Eq:AE-1} holds.
\end{lemma}

\begin{proof}
Since 
\[
\max_{x \in [0,1]} \frac{x(1-x)}{2x+1}  = \frac{1}{2}(2 - \sqrt{3}) \quad (\text{which is achieved at $x = \frac{1}{2}(\sqrt{3}-1)$}),
\]
and $7 - 3x = \frac{19}{4} + 3\big(\frac{3}{4} - x\big)$, 
it suffices to show 
\begin{align}
&\frac{26 - 9\sqrt{3}}{8} \big(\frac{3}{4} - x \big)^{3/2}  + \frac{57(2 - \sqrt{3})}{32 }\big(\frac{3}{4} - x \big)^{1/2}   \nonumber\\
		&\quad -   \frac{8 (1-x)^{1/2} (3x + 1)^{1/2}}{7x^{1/2}(2x+1)}\big[ (1 - x) - \frac{49}{32}\big(\frac{3}{4} - x \big)^3	
		\big] < 0 \text{ for } 0 < x \leq \frac{3}{4}.
		 \label{Eq:AE-2}
\end{align}

Let 
\begin{align*}
y(x) &= \big(\frac{3}{4} - x \big)^{1/2},\\
c(x) &= \frac{8(1-x)^{1/2} (3x + 1)^{1/2}}{7x^{1/2}(2x+1)},\\
H(y)
	&= \frac{4(26 - 9\sqrt{3}) y^3 +  57(2 - \sqrt{3}) y  }{8 + 32y^2 -  49y^6} \quad \text{ for } y \in [0,\sqrt{3}/2].
\end{align*}
Then \eqref{Eq:AE-2} is equivalent to 
\begin{equation}
H(y(x)) <   c(x) \text{ for } 0 < x \leq \frac{3}{4}.
	\label{Eq:AE-3}
\end{equation}
We make three observations.
\begin{enumerate}[(i)]
\item The function $y$ is strictly decreasing in $[0,3/4]$.
\item The function $c$ is strictly decreasing in $(0,1]$, as
\[
\frac{c'}{c} = -\frac{1 + 6x + 11 x^2 - 6x^3}{2x(1-x)(2x+1)(3x+1)} < 0.
\]

\item The function $H$ is strictly increasing in $[0,\sqrt{3}/2]$. To see this, factorize
\[
H = H_1 H_2 \text{ with } H_1(y) = \frac{y}{8 + 32y^2 -  49y^6} \text{ and } H_2(y) = 4(26 - 9\sqrt{3}) y^2 +  57(2 - \sqrt{3}).
\]
It is clear that $H_1$ and $H_2$ are positive. It is also clear that $H_2' \geq 0$ where equality holds only at $y = 0$. For $H_1$, we compute
\[
H_1'(y) = \frac{8 - 32 y^2 + 245 y^6}{(8 + 32 y^2 - 49 y^6)^2}.
\] 
By the arithmetic-geometric mean inequality, the numerator is positive:
\[
4 + 4 + 245 y^6 > 3 \times 3920^{1/3} y^2 \geq 32 y^2.
\]
Thus $H_1' > 0$. Hence $H$ is strictly increasing.
\end{enumerate}

\begin{table}[h]
\begin{center}
\caption{Arrangement to verify \eqref{Eq:AE-3}}
\label{Table1}
\begin{tabular}{c|c|c}
\hline
\hline
$[q_{j},q_{j+1}]$ & upper bound of $H(y(q_{j}))$ & lower bound of $c(q_{j+1})$ \\
\hline
$[0,.086]$ & $3.5555$ & $3.5655$\\
$[.086,.166]$ & $2.3472$ & $2.3538$\\
$[.166,.242]$ & $1.7875$ & $1.7906$\\
$[.242,.311]$ & $1.4561$ & $1.4580$\\
$[.311,.370]$ & $1.2420$ & $1.2449$\\
$[.370,.419]$ & $1.0973$ & $1.1000$\\
$[.419,.460]$ & $.9943$ & $.9949$\\
$[.460,.494]$ & $.9157$ & $.9166$\\
$[.494,.523]$ & $.8538$ & $.8550$\\
$[.523,.549]$ & $.8020$ & $.8032$\\
$[.549,.574]$ & $.7554$ & $.7562$\\
$[.574,.600]$ & $.7096$ & $.7097$\\
$[.600,.629]$ & $.6597$ & $.6604$\\
$[.629,.666]$ & $.5996$ & $.6009$\\
$[.666,.725]$ & $.5105$ & $.5119$\\
$[.725,.750]$ & $.2932$ & $.4758$
\end{tabular}
\end{center}
\end{table} 

In view of the above monotonicity of $y$, $c$ and $H$, to verify \eqref{Eq:AE-3}, we construct a finite increasing sequence $\{q_j\}_{j=1}^J$ with $q_1 = 0$ and $q_J = 3/4$ and check that
\[
H(y(q_{j})) < c(q_{j+1}) \text{ for } 1 \leq j \leq J-1.
\]
This sequence is constructed iteratively by choosing a rational $q_{j+1}$ to approximate $c^{-1}(H(y(q_j)))$, with suitable rounding to ensure the strict inequality $H(y(q_{j})) < c(q_{j+1})$. We summarize our arrangement in Table \ref{Table1} and skip the verification of the rational bounds for $H(y(q_{j}))$ and $c(q_{j+1})$ stated in the table.
\end{proof}

\begin{lemma}\label{Lem:Alg2}
It holds that
\begin{equation}
t(9 t^2 - 42t  + 17) < \frac{7}{4}  t^{-1/2} (1-t)^{1/2} (3t + 1)^{3/2}   \big(\frac{3}{4} - t\big)^{3/2}  \text{ for } t \in (0, 1/2).
	\label{Eq:AEX-1}
\end{equation}
\end{lemma}

\begin{proof}
We recast \eqref{Eq:AEX-1} as
\begin{multline}
h_1(t) := 9 t^2 - 42t  + 17  < \frac{7}{4}  t^{-3/2} (1-t)^{1/2} (3t + 1)^{3/2}   \big(\frac{3}{4} - t\big)^{3/2} =: h_2(t)\\
	  \text{ for } t \in (0, 1/2).
\label{Eq:AEX-2}
\end{multline}
We use the procedure in the proof of Lemma \ref{Lem:Alg1}. We first observe that $h_1$ and $h_2$ are both decreasing in $(0, 1/2)$: 
\begin{align*}
h_1'
	&= 18t - 42 < 0,\\
\frac{h_2'}{h_2} 
	&= -\frac{9 - 6 t + 41 t^2 - 48 t^3}{2 (-1 + t) t (1 + 3 t) (-3 + 4 t)} = -\frac{6 + 3(1-2t) + 17 t^2 + 24t^2(1 - 2t)}{2 (-1 + t) t (1 + 3 t) (-3 + 4 t)} < 0.
\end{align*}

\begin{table}[h]
\begin{center}
\caption{Arrangement to verify \eqref{Eq:AEX-2}}
\label{Table2}
\begin{tabular}{c|c|c}
\hline
\hline
$[t_{j},t_{j+1}]$ & upper bound of $h_1(t_j)$ & lower bound of $h_2(t_{j+1})$ \\
\hline
$[0,.180]$ & $17.0000$ & $17.0659$\\
$[.180,.252]$ & $9.7317$ & $9.7840$\\
$[.252,.301]$ & $6.9876$ & $6.9976$\\
$[.301,.346]$ & $5.1735$ & $5.1950$\\
$[.346,.403]$ & $3.5455$ & $3.5469$\\
$[.403,.500]$ & $1.5357$ & $1.7294$
\end{tabular}
\end{center}
\end{table}

Therefore, to prove \eqref{Eq:AEX-2}, we only need to construct a finite increasing sequence $\{t_j\}_{j=1}^J$ with $t_1 = 0$ and $t_J = 1/2$ such that
\[
h_1(t_j) < h_2(t_{j+1}) \text{ for } j = 1, \ldots, J-1.
\]
We iteratively choose a rational approximation $t_{j+1}$ to $h_2^{-1}(h_1(t_j))$, with suitable rounding to take care of the strict inequality $h_1(t_j) < h_2(t_{j+1})$. We summarize our arrangement in Table \ref{Table2} and skip the details.
\end{proof}
%----------------%

\section*{Rights retention statement.} 

For the purpose of Open Access, the authors have applied a CC BY public copyright licence to any Author Accepted Manuscript (AAM) version arising from this submission.

\section*{AI statement.}

ChatGPT-5.6 Sol was used to check for typographical errors. 

\section*{Acknowledgment.}

The authors would like to thank Xavier Lamy for his careful reading of the manuscript and for his useful comments. R.I. is partially supported by the ANR projects ANR-21-CE40-0004 and ANR-22-CE40-0006-01.
 
\def\cprime{$'$}

\end{document}